\documentclass[12pt]{article}
 \usepackage{amssymb,amsmath}
\def\vint_#1{-\kern-11pt\int_{#1}}
\begin{document}

\title{Regularity of weak solutions to the model Venttsel problem
for solutions of linear parabolic systems with nonsmooth in
 time principal matrix. $A(t)$-caloric method}


\author{Arina A. Arkhipova }
\medskip


\maketitle{We consider a model Venttsel type problem for linear
parabolic systems of equations. The Venttsel type boundary
condition  is fixed on the flat part of the lateral surface of a
given cylinder. It is defined by parabolic operator (with respect
to the tangential derivatives) and the conormal derivative.  The
H\"{o}lder continuity of a weak solution of the problem is proved
under optimal assumptions on the data. In particular, only
boundedness in the time variable of the principal matrices of the
system and the boundary operator is assumed.  All results are
obtained by so-called A(t)-caloric method \cite{AJS}.

 \section{Introduction}

In this paper we study regularity of weak solutions of the linear
parabolic systems under the Venttsel boundary condition on the
flat  part of  the lateral  surface of a given cylinder.

Let
 $B_1(0)=\{x\in \mathbb{R}^n:\,|x|<1\},$   $B^+_1(0)=B_1(0)\cap\{x_n>0\}$,
   $Q_1(0)=B_1(0)\times (-1,0)$,  and  $Q_1^+(0)=Q_1(0)\cap
 \{x_n>0\}$.

 We consider  a solution  $u:Q^+_1(0)\to  \mathbb{R}^N,\,\,N>1,$
 of the problem

 \begin{equation}\label{1.1}
u_t - \hbox{div}(a(z)\nabla u) = f(z),\,\,\,\, z \in Q^+_1(0),
\end{equation}
\begin{equation}\label{1.2}
u_t- \hbox{div}' (b(z')\nabla' u) +\frac{\partial
u}{\partial\mathbf{n}_a}=\psi(z'),\,\,\,\,z'\in \Gamma_1(0),
\end{equation}
 where  $\Gamma_1(0)=Q_1(0)\cap\{x_n=0\}$,
 $x=(x',x_n),\,z'=(x',0,t),$\,\,$\hbox{div}'(\cdot)=\Sigma_{i=1}^{n-1}\frac{\partial(\cdot)}{\partial
 x_i}$,  $\nabla'u=(u_{x_1},...,u_{x_{n-1}}).$

We  assume  that    $a(z)$  and  $b(z')$  are bounded positive
definite nondiagonal  $[nN\times nN]$  and $[(n-1)N\times(n-1)N]$
matrices,  $\frac{\partial u}{\partial \mathbf{n}_a}= (a(z)\nabla
u,\mathbf{n})$- vector of the conormal derivative.    Let  $f\in
L^2(Q_1^+(0))$,   $\psi\in L^2(\Gamma_1(0))$,  more exact
assumptions  on the data  will be formulated  below in conditions
{\bf H1} - {\bf H5}.
\medskip

\textbf{Definition\,1.1\,}\textit{ A function  $u\in H:=
L^2(\Lambda_1;W^1_2(B^+_1))\cap L^2(\Lambda_1; W^1_2(\gamma_1))$
is a weak solution to  problem  (\ref{1.1}),  (\ref{1.2})  if  it
satisfies  the identity
\begin{equation}\label{1.3}
\int_{Q^+_1}[-u\,\eta_t+(a(z)\nabla u,\nabla
\eta)]\,dz+\int_{\Gamma_1}[-u\,\eta_t+b(z')\nabla'u,\nabla'\eta]\,d\Gamma
=\int_{Q^+_1}f\,\eta\,dz+\int_{\Gamma_1}\psi\,\eta\,d\Gamma,
\end{equation}
with any  function} $\eta\in
C^\infty(\overline{Q^+_1}),\,\,\hbox{supp} \, \eta\in Q^+_1\cup
\Gamma_1.$

Here and below  we denote $\Lambda_1=(-1,0), \gamma_1=B_1(0)\cap
\{x_n=0\}$.  Certainly,  we  can  assume that  the test-functions
$\eta$  belong to the  space  $W^1_2(Q^+_1),\,\,\eta|_{\partial
Q^+_1\setminus \Gamma_1}= 0 $.

\medskip

The boundary condition (\ref{1.2}) includes both  the conormal
derivative  of  $u$ and parabolic  second order operator with  an
elliptic operator relative to the tangential derivatives. Such
boundary condition is referred to the Venttsel condition.  A
specific of the problem under consideration is that the boundary
operator is very strong one, and integral identity (\ref{1.3}) is
not homogeneous one with respect to similarity transformation of
the variables.

\medskip

To our days the \textit{scalar} ($N=1$) Venttsel problem
(considered in \cite{Ven} under more general boundary condition)
was studied  for  elliptic  and  parabolic nonlinear operators of
different classes  (see \cite{Naz}  and references therein).

In the case of vector-functions  ($N>1$), regularity of weak
solutions to the stationary Venttsel problem  was  studied  in
\cite{ALuk1}  and  \cite{ALuk2}. Elliptic operators with constant
coefficients  were considered in  \cite{ALuk1} and the Campanato
integral estimates were obtained for smooth solutions of the
Venttsel problem. Regularity of weak solutions of the
\textit{linear}  elliptic Venttsel problem  was studied in the
scale  of  the Morrey- Campanato spaces  in  \cite{ALuk2}. In
particular, the  H\"{o}lder continuity of  solutions  and their
first and second  derivatives  were proved in \cite{ALuk2},\,
$N>1$.

In  \cite{Ark},  the author considered the Venttsel
 problem  for  \textit{quasilinear} elliptic operators,  $N>1$.
 It is well known that  one can expect only  \textit{partial} regularity of weak
 solutions of systems with quasilinear operators (\cite{GG}, \cite{Cam81}, \cite{Cam84},  \cite{GM}  and references therein).  Moreover,  it was proved
 that singularities may  be  concentrated  near the boundary even  under  trivial Dirichlet condition \cite{G}.
In \cite{Ark}, the author proved  partial regularity of weak
solutions  of the Venttsel problem  for  quasilinear elliptic
operators with matrices  $a(x,u)$  and $b(x',u)$ $(f,\psi=0)$  and
estimated the Hausdorff measure of admissible singularities in the
vicinity of the boundary. To study regularity, it was applied the
so-called A-harmonic method \cite{DuG}. The method allowed to
relax continuity assumptions on $a(x,u)$  and $b(x',u)$ in
variables $x$ and  and $x'$ respectively  to  the integral (VMO)
continuity conditions.
\medskip

Here we study regularity of weak solutions  (\ref{1.3})  of the
Venttsel problem
 for  linear parabolic operators  in the model setting
 (\ref{1.1}),  (\ref{1.2}) and apply so-called $A(t)$- caloric method.

The method of $A$-caloric approximation was proposed by F.Duzaar,
G. Mingione   in  \cite{DuzMi} to study regularity
  of  weak solutions to a wide class of nonlinear parabolic
 systems  (see also \cite{BoDuMi}). According to this
 method, one can estimate locally the difference  in $L^2$- norm between
 a smooth solution of the simplest parabolic system with the
 constant matrix $A$ and a solution of the nonlinear system under
 consideration.  The main $A$-caloric lemma was later modified in
 \cite{AJS}  to  $A(t)$-caloric lemma. Such modification allowed
 to compare  well enough  a solution of quasilinear parabolic system  with the principal non smooth in time
  matrix and  a solution of the model system with the principal
 matrix $A(t)$ where $A(t)$ is only bounded in t. The application of
 $A(t)$-caloric lemma  allowed to prove partial regularity of weak solutions
 of quasilinear parabolic systems  with
 matrices $a(x,t,u)$  which are VMO-smooth in x  and only bounded
 in time variable t. Then it was proved partial regularity of a solution of the
 Cauchy-Dirichlet problem for quasilinear parabolic systems under
 the same assumptions on the principal matrix \cite{AS}.

 It should be mentioned works \cite{Kr1} - \cite{DK3}  where
 different classes of nonlinear parabolic \textit{scalar} equations  and
 linear parabolic systems were studied with non smooth in time
 principal matrices.
\medskip

In this paper we formulate $A(t)$-caloric lemmas in an appropriate
form (Section 3) and prove H\"{o}lder continuity of weak solutions
defined in (\ref{1.3}) under relax conditions on the matrices
$a(x,t)$ and $b(x',t)$.  More exactly, we assume  only integral
continuity of these matrices in the space variables and
boundedness in the time variable t.

\medskip
\textbf{} We introduce the following notation

\medskip

$$B_R(x^0)=\{x\in\mathbb{R}^n:\,\,|x-x^0|<R\},\quad
B^+_R(x^0)=B_R(x^0)\cap\{x_n>x_n^0\}\quad S_R(x^0)=\partial
B_R(x^0),$$ $$ S^+_R(x^0)=S_R(x^0)\cap \{x_n>x_n^0\},\quad
\gamma_R(x^0)=B_R(x^0)\cap\{x_n=x_n^0\},\quad
\Lambda_R(t^o)=(t^0-R^2,t^0), $$
$$Q_R(z^0)=B_R(x^0)\times(\Lambda_R(t^0)),\quad
Q_R^+(z^0)=B^+_R(x^0)\times\Lambda_R(t^0),\quad
\Gamma_R(z^0)=\gamma_R(x^0)\times\Lambda_R(t^0);$$
\medskip

we write $Q_R,\,Q^+_R,\, \Gamma_R, \,B^+_R,\, S_R^+$\, if
$z^0=0$;\,\,$|D|_m=mes_m\,D$ - is the Lebesgue measure of $D$  in
$\mathbb{R}^m$,\, $\omega_n=mes_n B_1$;
\medskip

$$
\vint_D f\,dz=\frac{1}{|D|_{n+1}}\int_D f\,dz,\quad
v_{r,z^0}=\vint_{Q^+_1\cap
Q_r(z^0)}v\,dz,\,\,\,v_{r,x^0}(t)=\vint_{B^+_1\cap
B_r(x^0)}v(x,t)\,dx,\,\,\, v_{x_\alpha}=\frac{\partial v}{\partial
x_\alpha};$$
\medskip

 $\|v\|_{p,Q}$ is the norm in
$L^p(Q)$; $W^1_2(Q)$  is  the Sobolev space  of functions $v$ in
$L^2(Q)$ possessing weak derivatives $v_{x_i},\,v_t\in
L^2(Q)$,\,$i=1,...n$.

We denote by $L^{2,\lambda}(Q^+_1;\delta)$  and
 $\mathcal{L}^{2,\lambda}(Q^+_1;\delta)$  the Morrey and Campanato spaces
 in the
\textbf{parabolic} metric
$\delta(z^1,z^2)=\max\{|x^1-x^2|,|t^1-t^2|^{1/2}\}$
\cite{Cam66}\,:
$$
L^{2,\lambda}(Q^+_1;\delta)=\{v\in
L^2(Q^+_1):\,\,\|v\|_{2,\lambda;Q^+_1}<\infty\},\quad \lambda\in
(0,n+2],
$$
where
$$
\|v\|_{2,\lambda;Q^+_1}=\left(\sup_{\rho\leq 1,z^0\in
Q^+_1}\frac{1}{\rho^\lambda}\int_{Q_\rho(z^0)\cap
Q^+_1}|v|^2\,dz\right)^{1/2};
$$
the space $\mathcal{L}^{2,\lambda}(Q^+_1;\delta),\,\lambda\in
(0,n+4],$ is the space of functions from $L^2(Q^+_1)$ with the
finite norm
$$
\mathbf{I}v\mathbf{I}_{2,\lambda;Q^+_1}=\|v\|_{2,Q^+_1}+[v]_{2,\lambda;Q^+_1}
$$
where the seminorm
$$[v]_{2,\lambda;Q^+_1}=\sup_{\rho\leq 1,z^0\in
Q^+_1}\frac{1}{\rho^\lambda}\int_{Q^+_1\cap
Q_\rho(z^0)}|v(z)-v_{\rho,z^0}|^2\,dz<\infty.$$

\medskip
We put

$$
W^{1,0}(Q^+_R(z^0))=L^2(\Lambda_R(t^0);W^1_2(B^+_R(x^0))),\,\,
W^{1,0}(\Gamma_R(z^0))=L^2(\Lambda_R(t^0);W^1_2(\gamma_R(x^0)));
$$
$$H(Q^+_R(z^0))=W^{1,0}(Q^+_R(z^0))\cap W^{1,0}(\Gamma_R(z^0)), $$
and
$$
V^{(1)}(Q^+_R(z^0))=C(\overline{\Lambda_R(t^0)};L^2(B^+_R(x^0)))\cap
L^2(\Lambda_R(t^0);W^1_2(B^+_R(x^0))),$$
$$V^{(2)}(\Gamma_R(z^0))=C(\overline{\Lambda_R(t^0)};L^2(\gamma_R(x^0)))\cap
L^2(\Lambda_R(t^0);W^1_2(\gamma_R(x^0))),$$
$$V(Q^+_R(z^0))=V^{(1)}(Q^+_R(z^0))\cap V^{(2)}(\Gamma_R(z^0)).
$$
\medskip

We write  $v\in \mathbf{B}(Q)$  instead of  $v\in
\mathbf{B}(Q;\mathbb{R}^N))$ for the sake of brevity.  Different
constants depending on the data of the problem are denoted by
$c$,\, $c_i$.

\section{The main results}

\medskip

We  formulate the basic assumptions:
\medskip

\begin{itemize}
\item\, {\bf H1} \,The matrix  $a$  is defined in $Q^+_1(0)$  and
has measurable bounded entries. There are positive constants
$\nu,\mu$ such that
\begin{equation}\label{H1}
(a(z)\xi,\xi)= a_{ij}^{\alpha \beta}(z)\xi^i_{\alpha}
\xi^j_{\beta}\geq \nu |\xi|^2, \quad \xi \in \mathbb R^{nN},
\end{equation}
\begin{equation} \nonumber
|(a(z)p,q)|\leq \mu|p||q|, \quad p,q \in \mathbb R^{nN},
\end{equation}
for almost all $z\in Q^+_1$.

\item {\bf H2} \, $ a(\cdot,t)\in VMO(B^+_1)$  for almost all
$t\in \Lambda_1$  and
\begin{equation}\label{H2} \sup_{\rho\leq r,z^0\in
Q^+_1}\vint_{\Lambda_\rho(t^0)} \vint_{B_\rho(x^0)\cap
B^+_1}|a(y,t)-a_{\rho,x^0}(t)|^2\,dy dt=:q_a^2(r)\to 0,\,\,r\to 0,
\end{equation}
where
$$a_{\rho,x^0}(t)=\vint_{B^+_1\cap B_\rho(x^0)}a(y,t)\,dy.$$

 \item \,{\bf H3}\, The matrix $b$ is defined on $\Gamma_1$ and has
 measurable bounded entries.  There are positive constants $\nu_0$
 and $\mu_0$ such that
 \begin{equation}\label{H3}
(b(z')\xi,\xi)=b_{ij}^{\alpha\beta}(z')\xi_i^\alpha\xi_j^\beta\geq
\nu_0 |\xi|^2,\quad \xi\in \mathbb{R}^{(n-1)N}, \end{equation}
\begin{equation}\nonumber
|b(z')p,q)|\leq \mu_0 |\xi|^2,\quad p,q \in \mathbb{R}^{(n-1)N},
\end{equation}
for almost all $z'\in \Gamma_1$.

\item \,{\bf H4}\,  $b(\cdot,t) \in VMO(\gamma_1)$ for almost all
$t\in \Lambda_1$  and
\begin{equation}\label{H4} \sup_{\rho\leq r,\,z'\in
\Gamma_1}\vint_{\Lambda_\rho(t')} \vint_{\gamma_1\cap
\gamma_\rho(\xi')}|b(y',t)-b_{\rho,\xi'}(t)|^2\,d\,y'\,dt=:q_b^2(r)\to
0,\quad z'=(\xi',t'),\quad r\to 0. \end{equation}

\item  \,{\bf \,H5\,} The function $f\in
L^{2,\lambda}(Q^+_1)$,\,\,$\psi\in L^{2,\lambda}(\Gamma_1)$  where
$\lambda =n-3+2\alpha,$\,\,$\alpha\in(0,1),$\  and  $n\geq 3$.

\item  \,{\bf\,H5'} \, The function $f\in
L^{2,n-2+2\alpha}(Q^+_1;\delta),\,$  and  $\psi\in
L^{2,n-3+2\alpha}(\Gamma_1;\delta),$\,  $\alpha\in (0,1),\,\,n\geq
3.$
\bigskip
\end{itemize}

\bigskip

\textbf{Theorem 2.1} {\it Let assumptions {\bf H1-H5} hold, $n\geq
3,$ \, and $u\in H$ be a weak solution to problem  (\ref{1.1}),
(\ref{1.2}). Then}

\noindent 1)  $u^G\in C^\alpha(\overline{\Gamma_R};\delta)$,
$\nabla'u^G\in L^{2,n-1+2\alpha}(\Gamma_R;\delta)$   \it{ for any
$R<1$, \,\,$u^G(z')=u(x',0,t)$\, and}
\begin{align}
\label{E10'}
\|u^G\|^2_{\mathcal{L}^{n+1+2\alpha}(\Gamma_R(0);\delta)}+\|\nabla'u^G\|^2_{L^{2,n-1+2\alpha}(\Gamma_R(0);\delta)}\leq
c\,\mathbb{M}_0\end{align} where}
$$ \mathbb{M}_0=\|u\|^2_{W^{1,0}_2(Q^+_1(0))}+\|u^G\|^2_{W^{1,0}_2(\Gamma_1(0))}
+\|f\|^2_{L^{2,n-3+2\alpha}(Q^+_1;\delta)}+\|\psi\|^2_{L^{2,n-3+2\alpha}(\Gamma_1;\delta)}.$$

\bigskip

\textbf{Theorem 2.2}  {\it Let the assumptions {\bf H1-H4} and
{\bf H5'} hold,  $n\geq 3$.  Then additionally to the assertion of
Theorem 2.1 \, $u\in
C^\alpha(\overline{Q^+_R(0)};\delta),$\,\,$\nabla\,u\in
L^{2,n+2\alpha}(Q^+_R(0);\delta)$  for any $R<1$ and}
\begin{align}\label{E10*}
\|u\|_{C^\alpha(\overline{Q^+_R(0)};\delta)}+\|\nabla\,u\|_{L^{2,n+2\alpha}(Q^+_R(0);\delta)}\leq
\,c\,\mathbb{L}_0\end{align}  {\it where}

$$ \mathbb{L}_0=\|u\|^2_{W^{1,0}_2(Q^+_1(0))}+\|u^G\|^2_{W^{1,0}_2(\Gamma_1(0))}
+\|f\|^2_{L^{2,n-2+2\alpha}(Q^+_1;\delta)}+\|\psi\|^2_{L^{2,n-3+2\alpha}(\Gamma_1;\delta)}.$$

\bigskip

\textbf{Theorem 2.3} \,\,{\it Let $n=2$,\, the assumptions {\bf
H1-H4} hold  and $u$ be a weak solution to problem  (\ref{1.1}),
(\ref{1.2}).

 \noindent{\bf I.}\,\, Let $f\in
 L^{2,\lambda_0}(\Gamma_1;\delta),\,\,\psi\in
 L^{2,\lambda_0}(\Gamma_1;\delta)$ with some }$\lambda_0\in [0,1)$.

  {\it Then}
\medskip

  i)\,\, $u^G\in C^{1/2}(\Gamma_{1-q};\delta),\,\,\nabla'u^G\in
  L^{2,2}(\Gamma_{1-q};\delta)$ {\it if}   $\lambda_0=0$;

  ii)\,\, $u^G\in
  C^{\alpha_0}(\Gamma_{1-q};\delta)$,\,\,$\nabla'u^G\in
  L^{2,2+\lambda_0}(\Gamma_{1-q};\delta),\,$
  $\alpha_0=\frac{\lambda_0+1}{2}$  {\it  if } $\lambda_0\in
  (0,1)$.

  \medskip

  \noindent{\bf II.}{\it Let $f\in
  L^{2,\lambda_0+1}(Q^+_1;\delta),\,\,\psi\in
  L^{2,\lambda_0}(\Gamma_1;\delta),\,\,\lambda_0\in (0,1)$.

  Then}

  \medskip
$u\in C^{\alpha_0}(\overline{Q^+_{1-q}};\delta),\,\,\nabla\,u\in
L^{2,2+\lambda_0}(Q^+_{1-q};\delta),\,\,\alpha_0=\frac{\lambda_0+1}{2}$.

{\it In all the  assertions number $q\in (0,1)$ is fixed
arbitrarily.}

\bigskip

\section{Auxiliary results}

\bigskip

 In this section we introduce  versions of the Caccioppoli
and Poincar\'e inequalities for weak solutions of the problem and
formulate $A(t)$-caloric lemmas in an appropriate form.

\bigskip

{\bf Lemma 3.1}(The Caccioppoli's inequality)  {\it Assume that
the conditions {\bf H1}  and {\bf H3} hold, $f \in L^2(Q^+_1)$,
$\psi\in L^2(\Gamma_1)$. Then  for all $z^0\in
\Gamma_1,\,\,Q^+_R(z^0)\subset Q^+_1$,   $l\in \mathbb{R}^N,$  a
solution  $u$  of (\ref{1.3}) satisfies the inequality

\begin{eqnarray}\label{3.1}
\int_{Q^+_{R/2}(z^0)}|\nabla\,u|^2\,dz+\int_{\Gamma_{R/2}(z^0)}|\nabla'u|^2\,d\Gamma\leq
 \frac{c}{R^2}\int_{Q^+_R(z^0)}|u-l|^2\,dz+\int_{\Gamma_R(z^0)}|u-l|^2\,d\Gamma \notag \\+
 c\,R^2\int_{Q^+_R(z^0)}|f|^2\,dz+
 c\,R^2\int_{\Gamma_R(z^0)}|\psi|^2\,d\Gamma.
\end{eqnarray}
The constants $c$ depend on the parameters from conditions {\bf
H1} and  {\bf H3}.
 }
\bigskip

{\bf Proof}\,\, It is not difficult to check that any solution of
(\ref{1.3}) belongs the class $V$ (see the notation in Section 2).
We omit smoothing in time  the  Steklov average procedure and put
formaly in (\ref{1.3})  $\eta=(u-l)\xi^2(x)\theta(t)$ where
$\xi(x)$ is a cut-off function for $B_R(x^0)$,  $\xi=1$  in
$B_{R/2}(x^0)$, $\theta\in C^1_0(\Lambda_R(t^0)),$\,\,$\theta=1$
in $\Lambda_{R/2}(t^0)$,  \, $z^0\in \Gamma_1$.  Inequality
(\ref{3.1}) follows  by the standard way.
\bigskip

{\bf Remark 3.1} Let here and below $u^G(x',t):=u(x',0,t)$,
$u^0=u(z)-u^G(z')$ and
  $l=(u^G)_{R,z^0}$  in (\ref{3.1}).  Then
 it follows from (\ref{3.1}) that
 \begin{eqnarray}\label{3.2}
 \int_{Q^+_{R/2}(z^0)}|\nabla\,u|^2\,dz+\int_{\Gamma_{R/2}(z^0)}|\nabla'u|^2\,d\Gamma\leq
 \frac{c}{R^2}\int_{Q^+_R(z^0)}|u^0|^2\,dz+\frac{c}{R^2}\int_{\Gamma_R(z^0)}|u^G-u^G_{R,z^0}|^2\,d\Gamma
 \notag \\
  +c\,R^2\int_{Q^+_R(z^0)}|f|^2\,dz+
  c\,R^2\int_{\Gamma_R(z^0)}|\psi|^2\,d\Gamma.
  \end{eqnarray}

\bigskip
{\bf Lemma 3.2} ( The Poincar\'e inequality) {\it Let the
assumptions {\bf H1} and {\bf H3} hold,  $z^0\in \Gamma_1$,
$Q_R^+(z^0)\subset Q^+_1$. Then for a weak solution $u$ of problem
(\ref{1.1}), (\ref{1.2}) the following inequality hold:}
\begin{eqnarray}\label{3.3}
\int_{Q^+_{R/2}(z^0)}|u(z)-u_{R/2,z^0}|^2\,dz+\int_{\Gamma_{R/2}(z^0)}|u^G(z')-u^G_{R/2,z^0}|^2\,d\Gamma
\notag \\  \leq
 c\,R^2\left(\int_{Q^+_R(z^0)}|\nabla\,u|^2\,dz   +\int_{\Gamma_R(z^0)}|\nabla'u^G|^2\,d\Gamma\right)
  +c\,R^4\left(\int_{Q^+_R(z^0)}|f|^2\,dz+\int_{\Gamma_R(z^0)}|\psi|^2\,d\Gamma\right).\end{eqnarray}
{\it The constants $c$ in  (\ref{3.3}) depend on the parameters
from assumptions {\bf H1} and  {\bf H3}.}
\medskip

\textbf{Proof}. As was noted in Lemma 3.1, weak solution of
(\ref{1.3}) is a function from the class $V$. For a point $z^0\in
\Gamma_1$ and a cylinder $Q^+_R(z^0)$  we fix  $s\in
\Lambda_R(t^0)\setminus \Lambda_{R/2}(t^0)$  and $\tau\in(s,t^0)$.
We fix  the same cut-off function  $ \xi(x)$ as in the proof of
Lemma 3.1. Let $\chi_\varepsilon(t)$ be a piecewise linear
function, $\chi_\varepsilon=1$  for $t\in [s,\tau]$  and
$\chi_\varepsilon=0$ when  $t<s-\varepsilon$  and
$t>s+\varepsilon$.  To simplify the explanation, we  omit  the
Steklov average procedure and  put  in (\ref{1.3})  $\eta
=(u(z)-l)\xi^2(x)\chi_\varepsilon(t)$  with any constant  vector
$l\in \mathbb{R}^N$.  After simple calculations we turn
$\varepsilon$  to   zero, and obtain the inequality
\begin{eqnarray}\label{3.4}
\int_{B^+_R(x^0)}|u(x,t)-l|^2\xi^2\,dx\left|\,^{t=\tau}_{t=s}
+\int_{\gamma_R(x^0)}|u^G(x',t)-l|^2\xi^2\,d\gamma\right|^{t=\tau}_{t=s}+
\nu\,\int^\tau_s\int_{B_R(x^0)}|\nabla\,u|^2\xi^2\,dx\,dt \notag \\
+\nu_0\int^\tau_s\int_{\gamma_R(x^0)}|\nabla'u^G|\xi^2\,dx'\,dt
\leq \mu\,\int^\tau_s\int_{B_R(x^0)}|\nabla\,u|\xi
\,|u-l|\,|\nabla\,\xi|\,dx\,dt
   \notag \\ +\mu_0\,\int^\tau_s\int_{\gamma(x^0)}|\nabla'u^g|\xi\,|u^G-l|\,|\nabla'\xi|\,d\gamma\,dt
  +\int^\tau_s\int_{B_R(x^0)}|f|\,|u-l|\xi^2\,dx\,dt
 \notag \\ +\int^\tau_s\int_{\gamma_R(x^0)}|\psi|\,|u^G-l|\xi^2\,d\gamma\,dt.
\end{eqnarray}
Now we put   $l=u^G_{R,x^0}(s)$  and estimate integrals with the
constants $\mu$ and  $\mu_0$ in the way:
\begin{eqnarray*}
\mu\int^\tau_s\int_{B^+_R(x^0)}|\nabla\,u|\,\xi\,|u-u^G_{R,x^0}(s)|\,|\nabla\,\xi|\,dx\,dt\leq
\delta\,\sup_{t\in(s,t^0)}\int_{B^+_R(x^0)}|u(x,t)-u^G_{R,x^0}(s)|^2\xi^2\,dx
\notag \\
+ c/\delta\int^\tau_s\int_{B^+_R(x^0)}|\nabla\,u|^2\,dx\,dt,
\end{eqnarray*}
\begin{eqnarray*}
\mu_0\int^\tau_s\int_{\gamma_R(x^0)} |\nabla'u^G| \,\xi\,
 |u^G-u^G_{R,x^0}(s)| \,|\nabla'\xi|\,d\gamma\,dt \leq \delta\sup_{t\in (s,t^0)}
\int_{\gamma_R(x^0)}|u^G-u^G_{R,x^0}(s)|^2\,\xi^2\,d\gamma \notag
\\
+c/\delta\,\int^\tau_s\int_{\gamma_R(x^0)}|\nabla'u^G|^2\,d\gamma\,dt.
\end{eqnarray*}
Further,
$$\int^\tau_s\int_{B^+_R(x^0)}|f|\,|u-u^G_{R,x^0}(s)|\,\xi^2\,dx\,dt\leq\,\delta
\sup_{t\in
(s,\tau)}\int_{B^+_R(x^0)}|u-u^G_{R,x^0}(s)|^2\xi^2\,dx+
\frac{c}{\delta} \,R^2\,\int_{Q^+_R(z^0)}|f|^2\,dz,
$$
and the integral with the function $\psi$ is estimated in the same
way.

 At last, by the Friedrichs and Poincare inequalities
$$
\int_{B^+_R(x^0)}|u(x,s)-u^G_{R,x^0}(s)|^2\,dx\leq
2\,\int_{B^+_R(x^0)}|u^0(x,s)|^2\,dx+2\,\int_{B^+_R(x^0)}|u^G(x',s)-u^G_{R,x^0}|^2\,d\,dx
$$
$$\leq c\,R^2\int_{B^+_R(x^0)}|\nabla\,u^0(x,s)|^2\,dx+
c\,R^3\,\int_{\gamma_R(x^0)}|\nabla'u(x',s)|^2\,d\gamma.$$

 Now we fix $\delta=1/8$ and derive from (\ref{3.4}) that
 \begin{eqnarray}\label{3.5}
\int_{B^+_R(x^0)}|u(x,\tau)- u^G_{R,x^0}(s)|^2\xi^2\,dx +
\int_{\gamma_R(x^0)}|u^G(x',\tau)-u^G_{R,x^0}(s)|^2\xi^2\,d\gamma
  \notag \\ \leq 1/2\,\sup_{t\in(s,t^0)}\int_{B^+_R(x^0)}|u(x,t)-u^G_{R,x^0}(s)|^2\xi^2\,dx  +1/2
\sup_{(t\in
(s,t^0)}\int_{\gamma_R(x^0)}|u^G(x',t)-u^G_{R,x^0}(s)|^2\,d\gamma
 \notag \\ + c\,R^2\int_{B^+_R(x^0)}|\nabla\,u^0|^2\,dx+
c\,R^3\,\int_{\gamma_R(x^0)}|\nabla'u^G(x',s)|^2\,d\gamma+
c\int_{Q^+_R(z^0)}|\nabla\,u|^2\,dz  \notag \\+
c\,\int_{\Gamma_R(z^0)}|\nabla'u^G(x's)|^2\,d\Gamma+
c\,R^2\int_{Q^+_R(z^0)}|f|^2\,dz+c\,R^2\int_{\Gamma(z^0)}|\psi|^2\,d\Gamma.
\end{eqnarray}
Taking  supremum in $\tau\in (s,t^0)$  in the left hand side of
(\ref{3.5}) we  obtain the inequality
\begin{eqnarray}\label{3.6}
\sup_{\tau\in
(s,t^0)}\int_{B^+_R(x^0)}|u(x,\tau)-u^G_{R,x^0}(s)|^2\,dx+
\sup_{\tau\in
(s,t^0)}\int_{\gamma_R(x^0)}|u^G(x',\tau)-u^G_{R,x^0}(s)|^2\,d\gamma
  \notag \\
  \leq c\,R^2\int_{B^+_R(x^0)}|\nabla\,u^0(x,s)|^2\,dx +
  c\,R^3\,\int_{\gamma_R(x^0)}|\nabla'u(x',s|^2\,d\gamma +
  c\int_{Q^+_R(z^0)}|\nabla\,u|^2\,dz \notag \\
 + c\,\int_{\Gamma_R(z^0)}|\nabla'\,u^G|^2\,d\Gamma +
  c\,R^2\int_{Q^+_R(z^0)}|f|^2\,dz  +
  c\,R^2\int_{\Gamma_R(z^0)}|\psi|^2\,d\Gamma.
  \end{eqnarray}
  To estimate the left hand side of (\ref{3.6}),  we use the fact
  that any  function  $\Phi(c)=\int_D|v-c|^2\,d\,D $  takes its
  minimum in $c=(v)_D$.  Then the inequality follows:
  \begin{eqnarray}\label{3.7}
  \sup_{t\in\Lambda_{R/2}(t^0)}\int_{B^+_{R/2}(x^0)}|u(x,t)-u_{R/2,x^0}(t)|^2\,dx
  + \sup_{t\in
  \Lambda_{R/2}(t^0)}\int_{\gamma_{R/2}(x^0)}|u^G(x',t)-u^G_{R/2,x^0}(t)|^2\,d\gamma
\notag \\ \leq c\,R^2\int_{B^+_R(x^0)}|\nabla\,u^0(x,s)|^2\,dx+
c\,R^3\int_{\gamma_R(x^0)}|\nabla'u^G(x',s)|^2\,d\gamma+ c
\int_{Q^+_R(z^0)}|\nabla\,u|^2\,dz \notag \\
+c\,R^2\int_{\Gamma_R(z^0)}|\nabla'u^G(z')|^2\,d\Gamma+
 c\,R^2\int_{Q^+_R(z^0)}|f|^2\,dz+
 c\,R^2\int_{\Gamma_R(z^0)}|\psi|^2\,\Gamma.
 \end{eqnarray}
The last inequaity we integrate over the interval
$\Lambda_R(t^0)\setminus \Lambda_{R/2}(t^0)$ and divide by the
measure of this interval. In a result we have
\begin{eqnarray}\label{3.8}
 \sup_{t\in
\Lambda_{R/2}(t^0)}\int_{B^+_{R/2}(x^0)}|u(x,t)-u_{R/2,x^0}(t)|^2\,dx+
\sup_{\Lambda_{R/2}(t^0)}\int_{\gamma_{R/2}(x^0)}|u^G(x',t)-u^G_{R/2,x^0}(t)|^2\,d\gamma
\notag \\
\leq
c\,\left(\int_{Q^+_R(z^0)}|\nabla\,u|^2\,dz+\int_{\Gamma_R(z^0)}|\nabla'u^G|^2\,d\Gamma
\right) + c\,R^2\left(\int_{Q^+_R(z^0)}|f|^2\,dz+
\int_{\Gamma_R(z^0)}|\psi|^2\,d\Gamma\right).
\end{eqnarray}
 Note that for a fixed $  s \in
\Lambda_R(t^0)\setminus\Lambda_{R/2}(t^0)$

$$ \int_{Q_{R/2}(z^0)}|u(z)-u_{R/2,z^0}|^2\,dz\leq
\int_{Q^+_{R/2}(z^0)}|u(z)- u_{R,x^0}(s)|^2\,dx\leq
\frac{R^2}{4}\sup_{t\in\Lambda_{R/2}(t^0)}\int_{B^+_{R/2}(x^0)}
|u(x,t)-u_{R,x^0}(s)|^2\,dx $$ and integral

$$\int_{\Gamma_{R/2}(z^0)}|u^G(z')-u^G_{R/2,z^0}|^2\,d\Gamma$$
 we estimate in the same way.

Now estimate (\ref{3.3})  follows  from (\ref{3.6}).
\bigskip

As was said in the Introduction, we apply in this work
$A(t)$-caloric method. Here we introduce two assertions in an
appropriate form  \cite{AJS},  \cite{AS}.\medskip

\textbf{Lemma 3.3}  {\it  Let positive numbers $\nu_0<\mu_0$ be
fixed. Suppose that  $m=(n-1)\,N$  and  $A(t)$   is  $[m\times m]$
matrix,  $A\in L^\infty (\Lambda_R(z^0))$,  satisfying   for
almost all $t\in \Lambda_R(t^0)$ the conditions
$$  (A(t)\xi,\xi)\geq \nu_0\,|\xi|^2,\quad
\quad |(A(t)\xi,\eta)|\leq \mu_0\,|\xi|\,|\eta|,\,\,\,\forall
\xi,\eta\in \mathbb{R}^m.
$$
Let a function $u^\Gamma\in W^{1,0}(\Gamma_R(z^0))$ be fixed.
 For any $\varepsilon>0$
there exist  a constant $C_\varepsilon=C(\varepsilon, \nu_0,\mu_0,
n, N)>0$, an $A(t)$-caloric function $h^G$ in
$W^{1,0}(\Gamma_{R/2}(z^0))$, and a function $\phi\in
C^1_0(\Gamma_R(z^0)),$ $\sup_{\Gamma_R(z^0)}|\nabla'\phi|\leq 1$
such that
\begin{eqnarray}\label{3.9}
\vint_{\Gamma_{R/2}(z^0)}(|h^G(z')-h^G_{R/2,z^0}|^2 +
R^2\,|\nabla'h^G|^2)\,d\Gamma\leq
2^{n+2}\,\vint_{\Gamma_R(z^0)}(|u^G(z')-u^G_{R,z^0}|^2+R^2\,|\nabla'u^G|^2)\,d\Gamma;
\end{eqnarray}
\begin{eqnarray}\label{3.10}
\vint_{\Gamma_{R/2}(z^0)}|u^G-h^G|^2\,d\Gamma\leq
\varepsilon\left(\vint_{\Gamma_R(z^0)}(|u^G-u^G_{R,z^0}|^2+R^2|\nabla'u^G|^2)\,d\Gamma\right)+
C_\varepsilon\,R^2\,\mathcal{L}^2_b(R,\phi)
\end{eqnarray}

 where}
\begin{eqnarray}\label{3.11}
\mathcal{L}_b(R,\phi)=\left|\vint_{\Gamma_R(z^0)}(-u^G\,\phi_t+
(A(t)\nabla'u^G,\nabla'\phi))\,d\Gamma\right|. \end{eqnarray}

\bigskip

\textbf{Lemma 3.4}
 {\it Let positive numbers $\nu<\mu $  be fixed,\,  $m=n\,N$.  Let
 an $[m \times m]$ matrix
 $A(t)$ satisfy the conditions of Lemma 3.3 with the parameters  $\nu,\,\mu$  and  $m=n\,N$.
 Suppose that a function $u\in W^{1,0}(Q^+_R(z^0))$,
 $u|_{\Gamma_R(z^0)}=0$.   For any $\varepsilon >0$  there exist
  a constant  $C_\varepsilon=C(\varepsilon, \nu,\mu,  m), \,$   an $A(t)$-caloric function  $h$ in $Q^+_{R/2}(z^0)$,
  \,$h|_{\Gamma_{R/2}(z^0)}=0$,  and a function  $\phi\in
  C^1_0(Q^+_R(z^0))$, \, $\sup_{Q^+_R}|\nabla\,\phi|\leq 1$  such  that
  \begin{eqnarray}\label{3.12}
  \vint_{Q^+_{R/2}(z^0)}(|h|^2+R^2\,|\nabla\,h|^2)\,dz\leq
  2^{n+2}\vint_{Q^+_R(z^0)}(|u|^2+R^2\,|\nabla\,u|^2)\,dz;
  \end{eqnarray}

\begin{eqnarray}\label{3.13}
\vint_{Q^+_{R/2}(z^0)}|u-h|^2\,dz\leq
\varepsilon\,\vint_{Q^+_R(z^0)}(|u|^2+R^2\,|\nabla\,u|^2)\,dz +
C_\varepsilon\,R^2\mathcal{L}^2_a(R,\phi),
\end{eqnarray}
where}
\begin{eqnarray}\label{3.14}
\mathcal{L}_a(R,\phi)=\left|\vint_{Q^+_R(z^0)}(-u\,\phi_t+
(A(t)\nabla\,u,\nabla\,\phi))\,dz\right|.
\end{eqnarray}

\bigskip

\textbf{Remark 3.2}\,\, It was proved in \cite{AS} that any
$A(t)$- caloric function  $h$  in
$Q^+_R(\xi)$,\,\,$h|_{\Gamma_R(\xi)}=0$,\,\,$\xi\in \Gamma_1$,
satisfies Campanato type  integral inequalities. We introduce here
the following estimate from Lemma 4\, \cite{AS}:
\begin{align}\label{Cam1}
\vint_{Q^+_\rho(\xi)}|h|^2\,dz\leq \left(
\frac{\rho}{r}\right)^2\vint_{Q^+_r(\xi)}|h|^2\,dz,\quad \rho\leq
r\leq R;\end{align}
 If $h^G$ is  $A(t)$-caloric function in
$\Gamma_R(z^0)$ then  (see \cite{AJS})
\begin{align}\label{Cam2}
\vint_{\Gamma_\rho(z^0)}|h-h_{\rho,z^0}|^2\,d\Gamma\,\leq
\,c\,\left(\frac{\rho}{r}\right)^2\,\vint_{Q^+_r(z^0)}
|h-h_{r,z^0}|^2\,d\Gamma,\quad \rho\leq r\leq R. \end{align}

\bigskip

\section{ H\"{o}lder continuity of $u$ on  $\Gamma_1$.  Proof of Theorem 2.1}

\bigskip

First, we apply  Lemma 3.3 to estimate  the function
$u^G(x',t)=u(x',0,t)\in W^{1,0}(\Gamma_R(z^0))$.
 Here $z^0=(x^0,t^0)\in \Gamma_1$ and  $\Gamma_{2R}(z^0)\subset \Gamma_1$  are fixed
 arbitrarily.

  We put $A^G(t)=\vint_{\gamma_R(x^0)}
 b(x',t)\,d\gamma$.  The matrix  $A^G(t)$ satisfies the conditions of Lemma 3.3 with
  the parameters  $ \nu_0<  \mu_0$.

  Now we fix an $\varepsilon>0$.  By Lemma 3.3, there exist a constant  $C_\varepsilon>0$, an  $A^G(t)$ caloric function  $h^G$  on $\Gamma_{R/2}(z^0)$,
  and  a function $\phi\in C^1_0(\Gamma_R(z^0)),$\,
  $\sup_{\Gamma_R(z^0)}|\nabla'\phi|\leq 1$   such that
  relations  (\ref{3.9})  and  (\ref{3.10})  are valid.
   We put $\eta(z)=\phi(z')\,m(x_n)$ in (\ref{1.3})\, where  $m\in C^1[0,R]$,  $m(0)=1,\,\,m(R)=0$.  Now  we
  estimate the expression  $\mathcal{L}_b(R,\phi)$  defined
  in  (\ref{3.11}):
            \begin{eqnarray}\label{4.1}
  \mathcal{L}_b(R,\phi) \leq
  \vint_{\Gamma_R(z^0)}|\Delta\,b|\,|\nabla'u^G|\,|\nabla'\phi|\,d\Gamma+
  R\,\vint_{Q^+_R(z^0)}|u\,\eta_t-
  \,(a(z)\nabla\,u,\nabla\,\eta)|\,dz  \notag \\ +
  R\,\vint_{Q^+_R(z^0)}|f\,\eta|\,dz
  +\vint_{\Gamma_R(z^0)}|\psi\,\phi|\,d\Gamma
  \end{eqnarray}
where
$$\Delta\,b= b_{R,x^0}(t)- b(z').$$
Taking into account that   $|\phi(z')|\leq c \,R$,  $|\phi_t|\leq
 c/R$, $|\eta_t|\leq c/R$,  we estimate the right hand side of (\ref{4.1})  and  derive the inequality

\begin{eqnarray}\label{4.2}
\mathcal{L}^2_b(R,\phi)\leq
c\vint_{Q^+_R(z^)}|u^0|^2\,dz+c\vint_{\Gamma_R(z^0)}|u^G-u^G_{R,z^0}|^2\,d\Gamma+
c\,R^4\vint_{Q^+_R(z^0)}|f|^2\,dz
+c\,R^2\vint_{\Gamma_R(z^0)}|\psi|^2\,d\Gamma \notag \\ +
c\,R^2\vint_{Q^+_R(z^0)}|\nabla\,u|^2\,dz +
c\,q^2_b(R)\vint_{\Gamma_R(z^0)}|\nabla'u^G|^2\,d\Gamma,\quad
u^0(z)=u(z)-u^G(z').
\end{eqnarray}
We estimate  two last terms in  relation  (\ref{4.2}) by
(\ref{3.2}). Then we multiply new relation by $R^2$ and obtain the
inequality
\begin{eqnarray}\label{4.3}
R^2\,\mathcal{L}^2_b(R,\phi)\leq c\,
q_b^2(R)\,R\,\vint_{Q^+_{2R}(z^0)}|u^0|^2\,dz +
c\,(q^2_b(R)+R)\vint_{\Gamma_{2R}(z^0)}|u^G-u^G_{2R,z^0}|^2\,d\Gamma
\notag \\ +
c\,R^5\vint_{Q^+_{2R}(z^0)}|f|^2\,dz+c\,R^4\vint_{\Gamma_{2R}(z^0)}|\psi|^2\,d\Gamma.
\end{eqnarray}
Now we introduce the function
\begin{eqnarray}\label{4.4}
\Phi(r,z^0)=r\,\vint_{Q^+_r(z^0)}|u^0|^2\,dz+\vint_{\Gamma_r(z^0)}|u^G(z')-u^G_{r,z^0}|^2\,d\Gamma,\quad
r\leq 2\,R.
\end{eqnarray}
Then  (\ref{4.3}) can be written in the form
\begin{eqnarray}\label{4.5}
R^2\,\mathcal{L}^2_b(R,\phi)\leq c\,(q_b^2(R)+R)\,\Phi(2R,z^0)+
c\, K(2R,z^0)
\end{eqnarray}
where  \begin{eqnarray}\label{kk}
 K(2R,z^0)=
R^5\,\vint_{Q_{2R}(z^0)}|f|^2\,dz+
R^4\,\vint_{\Gamma_{2R}(z^0)}|\psi|^2\,d\Gamma. \end{eqnarray}

 It follows from inequalities (\ref{3.9}), (\ref{3.10}), (\ref{4.5}) and (\ref{3.2})  that

\begin{eqnarray}\label{4.6}
\vint_{\Gamma_{R/2}(z^0}(|h^G-h^G_{R/2,z^0}|^2+R^2\,|\nabla'h^G|^2)\,d\Gamma\leq
c\,\Phi(2R,z^0)+c\,K(2R,z^0),
\end{eqnarray}
\begin{eqnarray}\label{4.7}
\vint_{\Gamma_{R/2}(z^0)}|u^G-h^G|^2\,d\Gamma\leq
\left(\varepsilon + C_\varepsilon(R+q_b^2(R))\right)\,\Phi(2R,z^0)
+C_\varepsilon\,c\, K(2R,z^0).
\end{eqnarray}

The next step of the proof is to apply  Lemma 3.4  to the function
$u^0$  in  $Q^+_R(z^0)$,   $u^0|_{\Gamma_R(z^0)}=0$. We put the
matrix $A(t)=\vint_{B^+_R(x^0)}a(x,t)\,dx$.  For the fixed earlier
$\varepsilon$,  there exist a constant $C^*_\varepsilon$, an
$A(t)$ caloric function  $h$ defined in $Q^+_{R/2}(z^0),\,$
 $h|_{\Gamma_{R/2}(z^0)}=0$,  and a function $\phi^*\in
 C^1_0(\overline{Q^+_R(z^0)}),\,\,\sup_{Q^+_R(z^0)}|\nabla\,\phi^*|\leq
 1$  such that
 inequalities  (\ref{3.12}), (\ref{3.13})  hold with  $u^0$
 instead of $u$.  Using identity (\ref{1.3}), we estimate the expression
 $\mathcal{L}_a(R,\phi^*)$ in the way
 \begin{eqnarray}\label{4.8}
 \mathcal{L}_a(R,\phi^*):=\left|\vint_{Q^+_R(z^0)}(-u^0\,\phi^*_t+(A(t)\nabla\,u^0,\nabla\,\phi^*))\,dz\right|
 \leq
 \vint_{Q^+_R(z^0)}|\Delta\,a||\nabla\,u^0|\,|\nabla\,\phi^*|\,dz
 \\ \notag
 + \vint_{Q^+_R(z^0)}|u^G\,\phi_t^* -
 (a(z)\,\nabla\,u^G,\nabla\,\phi^*)|\,dz
 + \vint_{Q^+_R(z^0)}|f\,\phi^*|\,dz.
 \end{eqnarray}
Recall that $|\phi^*|\leq c\,R$  and  $|\phi_t^*|\leq c/R. $ Thus,
$$
\left|\vint_{Q^+_R(z^0)}u^G\,\phi^*_t\,dz\right|=\left|\vint_{Q^+_R(z^0)}(u^G-u^G_{R,z^0})\,\phi^*_t\,dz\right|
\notag \\ \leq c\,
\,R^{-1}\vint_{\Gamma_R(z^0)}|u^G-u^G_{R,z^0}|\,d\Gamma,
$$
$$
\vint_{Q^+_R(z^0)}|f|\,|\phi^*|\,dz\leq
c\,R\vint_{Q^+_R(z^0)}|f|\,dz,\quad
\vint_{Q^+_R(z^0)}|a\,\nabla'u^G,\nabla\,\phi^*|\,dz\leq
c\,\vint_{\Gamma_R(z^0)}|\nabla'u^G|\,d\,\Gamma.
$$
Now the estimate follows
\begin{eqnarray}\label{4.9}
\mathcal{L}_a^2(R,\phi^*)\leq
\vint_{Q^+_R(z^0)}|\Delta\,a|^2\,dz\,\vint_{Q^+_R(z^0)}|\nabla\,u^0|^2\,dz
+ \frac{c}{R^2}\vint_{\Gamma_R(z^0)}|u^G- u^G_{R,z^0}|^2\,d\Gamma
  \notag \\ + c\,R^2\,\vint_{Q^+_R(z^0)}|f|^2\,dz
+\vint_{\Gamma_R(z^0)}|\nabla'u^G|^2\,d\Gamma.
\end{eqnarray}
Then we can apply condition {\bf H2} and inequality (\ref{3.2}) to
estimate the first and the last integrals   in the right hand side
of (\ref{4.9}).

In a result, we obtain the inequality
\begin{eqnarray}\label{4.10}
R^2\,\mathcal{L}^2_a(R,z^0)\leq c\,(q^2_a(R)\,R^{-1}
+1)[\Phi(2R,z^0) +\,K(2R,z^0)].
\end{eqnarray}
Now it follows from (\ref{3.13}) and (\ref{4.10})   that
\begin{eqnarray}\label{4.11'}
\vint_{Q^+_{R/2}(z^0)}|u^0-h|^2\,dz\leq \varepsilon
\vint_{Q^+_R(z^0)}(|u^0|^2+r^2\,|\nabla\,u^0|^2)\,dz +
C^*_\varepsilon\,R^2\,\mathcal{L}^2(R,\phi^*) \notag \\
 \leq c\,\{\varepsilon\, R^{-1}+ C^*_\varepsilon\,
 (q^2_a(R)\,R^{-1} + 1)\}\,\Phi(2R,z^0)+c\,C^*_\varepsilon\,(q^2_a(R)\,R^{-1}+1) K(2R,z^0).
 \end{eqnarray}
\medskip

On the next step we will use known estimates for $A(t)$ caloric
functions   $h^G$  and $h$  (see Remark 3.2)  and  estimates
(\ref{4.7}) and (\ref{4.11'}).

The following chain of the inequalities are valid  for $\rho\leq
R/2$: \begin{eqnarray}\label{4.11}
\Phi(\rho,z^0)=\rho\,\vint_{Q^+_\rho(z^0)}|u^0|^2\,dz+\vint_{\Gamma_\rho(z^0)}|u^G-u^G_{\rho,z^0}|^2\,d\Gamma\leq
2\,\rho\,\vint_{Q^+_\rho(z^0)}|u^0-h|^2\,dz \notag \\
+ \,2\,\rho\,\vint_{Q^+_\rho(z^0)}|h|^2\,dz
 +
 2\vint_{\Gamma_\rho(z^0)}|(u^G-h^G)-(u^G_{\rho,z^0}-h^G_{\rho,z^0})|^2\,d\Gamma
 +2 \vint_{\Gamma_\rho(z^0)}|h^G-h^G_{\rho,z^0}|^2\,d\Gamma.
 \end{eqnarray}
Applying estimates (\ref{Cam1})  and  (\ref{Cam2})  for $A(t)$-
caloric functions $h$ and $h^G$,  $\rho\leq R/2$,  we obtain from
(\ref{4.11}) the inequality
\begin{eqnarray}\label{4.12}
\Phi(\rho,z^0)\leq
\rho\,\left(\frac{R}{\rho}\right)^{n+2}\vint_{Q^+_{R/2}(z^0)}|u^0-h|^2\,dz
+c\,\left(\frac{R}{\rho}\right)^{n+1}\vint_{\Gamma_{R/2}(z^0)}|u^G-h^G|^2\,d\Gamma
\notag \\
 + \rho\,\left(\frac{\rho}{R}\right)^2\vint_{Q^+_{R/2}(z^0)}|h|^2\,dz
+
c\,\left(\frac{\rho}{R}\right)^2\vint_{\Gamma_{R/2}(z^0)}|h^G-h^G_{R/2,z^0}|^2\,d\Gamma.
\end{eqnarray}
Now we use   (\ref{3.12})  and  (\ref{3.2}) to obtain the
inequality
\begin{eqnarray}\label{4.13}
\vint_{Q^+_{R/2}(z^0)}|h|^2\,dz\leq c\,R^{-1} (\Phi(2R,z^0)+
K(2R,z^0)).
\end{eqnarray}
Further, applying relations  (\ref{4.6}), (\ref{4.7}),
(\ref{4.11'}) and (\ref{4.13}), we estimate the right-hand side of
(\ref{4.12}) with $r=2R$   as follows:
\begin{eqnarray}\label{4.14}
\Phi(\rho,z^0)\leq c_0\left\{\left(\frac{\rho}{r}\right)^2+
\varepsilon\left(\frac{r}{\rho}\right)^{n+1} + \hat{C}_\varepsilon
\left(\frac{r}{\rho}\right)^{n+1}[r + q_a^2(r) +
q_b^2(r)]\right\}\Phi(r,z^0)
\notag \\
 + \hat{C}_\varepsilon \,c\, K(r,z^0), \quad \rho\leq r/4, \quad  \hat{C}_\varepsilon= \max \{C^*_\varepsilon,\,C_\varepsilon\}.
 \end{eqnarray}
It follows from the assumptions  {\bf H5} on $f$, $\psi$ and the
definition (\ref{kk}) that
$$
K(r)\leq K_0\,r^{2\alpha},\quad  n\geq 3,$$
$$
K_0=\|f\|^2_{L^{2,\lambda}(Q^+_1;\delta)}+\|\psi\|^2_{L^{2,\lambda}(\Gamma_1;\delta)}.
$$

  We put now in (\ref{4.14})  $\rho=\tau\,r$
 with $\tau\leq 1/4$ to be chosen later.

 Then
 \begin{equation}\label{4.15}
 \Phi(\tau\,r,z^0)\leq c_0\{\tau^2+\varepsilon\,\tau^{-(n+1)}+
 \hat{C}_\varepsilon\,\tau^{-(n+1)}[r+q^2_a(r)+q_b^2(r)]\}\Phi(r,z^0)+c\,\hat{C}_\varepsilon\,K_0\,r^{2\alpha}.
 \end{equation}
Now we fix  $\beta=\frac{\alpha+1}{2},\,\beta>\alpha,$  and choose
$\tau$ to satisfy the relation \begin{equation}\label{4.16}
c_0\,\tau^2\leq \frac{\tau^{2\beta}}{3}.\end{equation} Further, we
fix $\varepsilon<1$  in the way
\begin{equation}\label{4.17}
c_0\,\varepsilon\,\tau^{-(n+1)}\leq
\frac{\tau^{2\beta}}{3}.\end{equation}
The parameters $\tau$ and
$\varepsilon$ are fixed by the data of the problem and do not
depend on $z^0$.

At last, we can specify the choice of $r_0=2R_0$  by requiring
that
\begin{equation}\label{4.18}
 c\,\hat{C}_\varepsilon\,\tau^{-(n+1)}[r_0+q^2_a(r_0)+q^2_b(r_0)]\leq
 \frac{\tau^{2\beta}}{3}.
 \end{equation}
In a result,
\begin{equation}\label{4.19}
\Phi(\tau\,r,z^0)\leq \tau^{2\beta}\Phi(r,z^0)+
c\,K_0r^{2\alpha}.\end{equation} Proceeding by induction  in
relation (\ref{4.19}) for $r_j=\tau^j\,r$,  $j\in \mathbb{N},$
 we
obtain the inequality
\begin{equation}\label{4.20}
\Phi(r^j,z^0)\leq \tau^{2\alpha\,j}(\Phi(r,z^0)+c_1\,K_0
\,r^{2\alpha}).
\end{equation}
We can assert now that
\begin{equation}
\Phi(\rho,z^0)\leq
c\,\left(\frac{\rho}{r}\right)^{2\alpha}\left[\Phi(r,z^0)+ K_0\,
r^{2\alpha}\right],\quad \forall \rho\leq r\leq r_0.
\end{equation}
Thus, \begin{eqnarray}\label{4.21}\sup_{\rho\leq
r_0}\rho^{-(n+1+2\alpha)}\left(\int_{Q^+_\rho(z^0)}|u^0|^2\,dz+\int_{\Gamma_\rho(z^0)}|u^G-u^G_{\rho,z^0}|^2\,d\Gamma\right)
 \notag \\ \leq
c(r_0^{-1})(\|\nabla\,u\|^2_{L^2(Q^+_1)}+\|\nabla'u^G\|^2_{L^2(\Gamma_1)}+
K_0)\leq c(r_0^{-1})\,\mathbb{M}_0,\end{eqnarray} where
$\mathbb{M}_0$ is defined in  (\ref{E10'}).
    Taking supremum in $z^0\in \Gamma_{1-q}(0)$ in the left-hand
    side of (\ref{4.21}), ($q\in (0,1)$ is any fixed number,\,  $r_0$ satisfies (\ref{4.18})  and  $r_0\leq q$), we
    obtain  the estimate of the seminorm of $u^G$ in
    $\mathcal{L}^{2,n+1+2\alpha}(\Gamma_{1-q}(0);\delta)$:
    \begin{eqnarray}\label{4.22}
    [u^G]_{\mathcal{L}^{2,n+1+2\alpha}(\Gamma_{1-q}(0);\delta)}\leq
    c(r_0^{-1},q^{-1})\,
    (\|\nabla\,u\|_{L^2(Q^+_1(0))}
    +\|\nabla'u^G\|^2_{L^2(\Gamma_1(0))} + K_0)\leq c\,\mathbb{M}_0.
    \end{eqnarray}
Due to the isomorphism between
$\mathcal{L}^{2,n+1+2\alpha}(\Gamma_{1-q}(0);\delta)$  and
 $C^\alpha(\overline{\Gamma_{1-q}(0)};\delta)$ in the parabolic metric, we
obtain estimate of $C^\alpha$-norm of $u^G$ in $\Gamma_{1-q}(0)$.
 Moreover, estimate (\ref{3.2}) provides that
$$
\|\nabla'u^G\|^2_{L^{2,n-1+2\alpha}(\Gamma_{1-q}(0);\delta)}\leq
c\,\mathbb{M}_0,
$$
and estimate (\ref{E10'}) follows. $\bullet$

\bigskip

\section{ Proofs  of Theorem 2.2  and 2.3}

\bigskip
Here we consider problem (\ref{1.1}), (\ref{1.2}) in the form
\begin{equation}\label{5.1}
u_t- div (a(z)\nabla\,u)= f(z),  \quad z\in Q^+_1,
\end{equation}
\begin{equation}\label{5.2}
u|_{\Gamma_1}=\phi(z), \end{equation}  where $ \phi(z)=u^G(x',t)$
 and $f\in L^{2,n-2+2\alpha}(Q^+_1;\delta)$. If all assumptions of Theorem 2.2  hold then   $ u^G(x',t)\in C^\alpha
(\Gamma_{1-q}(0);\delta),$ $ \forall q\in (0,1)$ by Theorem 2.1.
\medskip

Here we prove further smoothness results for $u$.
\medskip
First, we want to recall some known results on the regularity of
solutions (\ref{5.1}), (\ref{5.2}).
\medskip

\noindent {\bf Proposition  5.1.}  \textit{ Let $N=1$, \,$n\geq
2$, the assumptions {\bf H1--H4} hold,\, $f\in
L^p(Q^+_{R_0}),\,p>\frac{n+2}{2},\,\,n\geq 2,$  and  $\psi$
satisfies the condition {\bf H5}  with $\alpha=2-\frac{n+2}{p}>0$.
Then  $u\in C^\beta(\overline{Q_R^+};\delta)$ with some $\beta\leq
\alpha$ and $R<R_0$.}

\medskip

 This result is a
consequence of a weak form of the maximum principle (see, for
example  ~\cite{LSU}, Chapter III, \S 10). Indeed, if $f\in
L^p(Q^+_1)$  with $p>\frac{n+2}{2}$ then  it belongs to the space
$L^{2,n-2+2\alpha}(Q^+_1;\delta),$ $\alpha=2-\frac{n+2}{p}>0$.
Using Theorem 2.1,  we  can assert that $u^G\in
C^\alpha(\Gamma_{1-q}(0);\delta)$. Thus, by the mentioned integral
form of the maximum principle
 $u\in C^\beta(\overline{Q^+_{1-q}(0)};\delta)$,\,\,$\beta\leq \alpha$ in the
\textit{scalar} case.

\medskip

Let now $N>1$  and the assumptions of Theorem 2.2 hold. Analyzing
the proof of Theorem 1 in \cite{AJS} where regularity problem for
\textit{quasilinear} systems was studied, one can assert local
smoothness of weak solutions of systems (\ref{5.1}). More exactly,
the following proposition is valid.
\medskip

 \noindent{\bf Proposition  5.2.}\textit{Let the matrix  $a$  satisfy the conditions {\bf
 H1,\,H2},\,$n\geq 2$, \,\,
   $f\in L^{2,n-2+2\alpha}(Q^+_1;\delta)$,
  $\alpha\in (0,1)$,  and $u$ be a weak solution of (\ref{5.1}) from $H_1=W^{1,0}_2(Q^+_1)$.
   Then  $u\in C^\alpha(\overline{Q'};\delta),\,\,\forall Q'\subset\subset Q^+_1$,  and }
\begin{align}\label{5.3}
\|u\|_{C^\alpha(\overline{Q'};\delta)}+\|\nabla\,u\|_{L^{2,n+2\alpha}(Q';\delta)}\leq
c_1\,\|u\|_{H_1}+ c_2\,\|f\|_{L^{2,n-2+2\alpha}(Q^+_1;\delta)}.
\end{align}

 \textit{Moreover,  for any $\xi\in Q'$ and} $\rho\leq r\leq
 \delta(\xi,\partial _p(Q^+_1))$
 \begin{align}\label{5.4}
 \Phi(\rho,\xi):=\frac{1}{\rho^{n+2+2\alpha}}\int_{Q_\rho(\xi)}|u-u_{\rho,\xi}|^2\,dz\leq
 c_3\{\Phi(r,\xi)+\|f\|^2_{L^{2,n-2+2\alpha}(Q^+_1;\delta)}\}.
 \end{align}
 \textit{The constants $c_1  - c_3$ depend on the parameters from conditions {\bf H1,\,H2}, $\alpha$,
  and the constant $c_1$ also depends on $\delta(Q',\partial_pQ^+_1)>0$.}

\medskip

As a consequence of Proposition 5.2 and Theorem 2.1, we obtain the
following assertion.
\medskip

\noindent{\bf Proposition 5.3.}  \,\textit{Let the assumptions
{\bf H1--H4} hold.  If  $2\alpha-1>0$  in the assumption {\bf H5'}
then $u\in C^\beta(\overline{Q^+_{1-q}};\delta)$ with
$\beta=\alpha-1/2>0$  and any fixed $q\in (0,1)$.}

\medskip

Indeed,  if $z^0\in \Gamma_{1-q}(0)$ then by Theorem 2.1  for
$\rho\leq r_0\leq q$
$$\Psi(\rho,z^0):=
\frac{1}{\rho^{n+1+2\alpha}}\int_{Q_\rho(z^0)}|u-u_{\rho,z^0}|^2\,dz
 \leq\frac{1}{\rho^{n+1+2\alpha}}\int_{Q_\rho(z^0)}|u-u^G_{\rho,z^0}|^2\,dz
 $$
 \begin{align}\label{5.5}\leq
 \frac{2}{\rho^{n+1+2\alpha}}\int_{Q_\rho(z^0)}|u^0|^2\,dz+
 \frac{2\rho}{\rho^{n+1+2\alpha}}\int_{\Gamma_\rho(z^0)}|u^G-u^G_{\rho,z^0}|^2\,d\Gamma,
 \,\,\,u^0(z)=u(z)-u^G(z').
 \end{align}
The right hand side of (\ref{5.5})  can be estimated  by Theorem
2.1. Thus,
$$\Psi(\rho,z^0)\leq c\,\mathbb{L}_0$$
for any point $z^0\in \Gamma_{1-q}(0)$, here $r_0$ does not depend
on $z^0$  and $\mathbb{L}_0$ is defined by (\ref{E10*}).

The standard procedure of "sewing" together local inner and
boundary estimates  for $\Psi(\rho,\cdot)$  provides estimate of
this function for all $\xi\in \overline{Q_{1-q}^+(0)}.$  We remark
that $n+1+2\alpha=n+2+2\beta,\,$  $\beta=\alpha-1/2,\,$  and the
H\"{o}lder continuity  of $u$ in $\overline{Q^+_{1-q}(0)}$ follows
with the exponent $\beta=\alpha-1/2$.  We do not explain in
details the proof of Proposition 5.3  because below we prove the
more strong assertion of Theorem 2.2.
\bigskip

\noindent {\bf Proof of Theorem 2.2.}

\medskip

  We start with the transformation of problem (\ref{5.1}),
 (\ref{5.2}) to the homogeneous one.

 \medskip

  We put
 $u^0(z)=u(z)-u^G(z'),\,\,  u^0|_{\Gamma_1}=0,$ and  define a weak solution to the problem
\begin{align}\label{5.6}
u^0_t-div(a(z)\nabla\,u^0)=-u^G_t+ div(a(z)\nabla'u^G)+f(z),\quad
z\in Q^+_1(0),
\\ \notag u^0|_{\Gamma_1(0)}=0.
\end{align}

\noindent {\bf Definition 5.1.}\,  \textit{A function $u^0\in
W^{1,0}_2(Q^+_1(0)),\,\,u^0|_{\Gamma_1(0)}=0,$  is a weak solution
to  problem (\ref{5.6}) if it  satisfies the identity}
\begin{align}\label{5.7}
\int_{Q^+_1(0)}[-u^0\,\eta_t
+(a(z)\nabla\,u^0,\nabla\,\eta)]\,dz=\int_{Q^+_1(0)}[u^G\,\eta_t
-(a(z)\nabla'u^G,\nabla\,\eta)+f\,\eta]\,dz \end{align}
 for any  $\eta\in \overset{0}{W}^{1,1}_2(Q^+_1(0)).$

\medskip
To prove Theorem 2.2 it is enough to  state H\"{o}lder continuity
of $u^0$ in $\overline{Q^+_{1-q}(0)},\,\,q\in (0,1).$
\medskip

As a first step, we prove that there exist  half derivatives  in
$t$  of the functions $u$ and $u^G$.

\medskip

 We fix $z^0\in \Gamma_1(0)$ and $Q^+_{2R}(z^0)\subset
Q^+_1.$  Let $\omega(t)\in
C^1_0(\Lambda_{2R}(t^0)),$\,\,$\omega(t)=1$ in $\Lambda_{R}(t^0)$;
let  $d(x)$ be a cut-off function for $B_{2R}(x^0)$,  $d(x)=1$ in
$B_R(x^0).$  Note that $|\omega'(t)|\leq
\frac{c}{R^2},\,\,|\nabla\,d(x)|\leq \frac{c}{R}.$  We put
\begin{align}\label{5.8}
v(z)=(u(z)-u^G_{2R})\,\omega(t)\,d(x),\quad v|_{\partial
Q^+_{2R}(z^0)\setminus \Gamma_{2R}(z^0)}=0,
\end{align}
and prove the following proposition.

\noindent {\bf Proposition 5.4.} \,\,\textit{Let  assumptions
{\bf H1,\,H3} hold and $u$ be a weak solution to problem
(\ref{5.1}),\,(\ref{5.2})  in $Q^+_{2R}(z^0)\subset Q^+_1(0)$,
$z^0\in \Gamma_1(0)$.  Then}

\noindent 1)\, $v\in
H^{1/2}(\Lambda_{2R}(t^0);L^2(B^+_{2R}(x^0)),\,\,v^G=v|_{x_n=0}\in
H^{1/2}(\Lambda_{2R}(t^0);L^2(\gamma_{2R}(x^0)))$ \textit{ where
$v$ defined by (\ref{5.8}),  and}
\begin{align}\label{5.9}
\|v\|^2_{H^{1/2}(\Lambda_{2R}(t^0);L^2(B^+_{2R}(x^0)))}+\|v^G\|^2_{H^{1/2}(\Lambda_{2R}(t^0);L^2(\gamma_{2R}(x^0)))}
\\ \notag
\leq
c_1\,\{\|\nabla\,v\|^2_{2,Q^+_{2R}(z^0)}+\|\nabla'v^G\|^2_{2,\Gamma_{2R}(z^0)}
+R^{-2}(\|u^0\|^2_{2,Q^+_{2R}(z^0)}+\|u^G-u^G_{{2R},z^0)}|_{2,\Gamma_{2R}}^2)\\
\notag  +R^2(\|f\|^2_{2,Q^+_{2R}}+\|\psi\|^2_{2,\Gamma_{2R}})\};
\end{align}
\medskip

\noindent 2) \,\,$u\in H^{1/2}(\Lambda_R(t^0);L^2(B^+_R(x^0)))$,
\,$u^G\in H^{1/2}(\Lambda_R(t^0);L^2(\gamma_R(x^0)))$
\textit{and}
\begin{align}\label{5.10}
\|u\|_{H^{1/2}(\Lambda_R(t^0);L^2(B^+_R(x^0)))}^2+\|u^G\|^2_{H^{1/2}(\Lambda_R(t^0);L^2(\gamma_R(x^0)))}
\\ \notag
\leq
c_2\,\{\|\nabla\,u\|^2_{2,Q^+_{2R}}+\|\nabla'u^G\|^2_{2,\Gamma_{2R}}+R^{-2}(\|u^0\|^2_{2,Q^+_{2R}}
+\|u^G-u^G_{2R}\|_{2,\Gamma_{2R}}^2) \\ \notag
+R^2(\|f\|^2_{2,Q^+_{2R}}+\|\psi\|^2_{2,\Gamma_{2R}})\}.
\end{align}
\textit{ The constants $c_1$ and $c_2$ depend on the parameters
from the conditions {\bf H1,\,H3}, \, $n, \,N,\,$ and do not
depend on \,$z^0$ and \, $R$.}

 \medskip

\noindent{\bf Proof.} \,\,  We consider identity (\ref{5.3}) with
$\eta(z)=\omega(t)d(x)\xi(z)$ where $\omega(t)$ and $d(x)$ are the
same as in (\ref{5.8}). The function $\xi\in
W_2^{1,1}(\hat{Q})\cap W^{1,1}_2(\hat{\Gamma})$,  here and below
we denote
$$
\hat{Q}=B^+_{2R}(x^0)\times \mathbb{R}^1,\quad
\hat{\Gamma}=\gamma_{2R}(x^0)\times \mathbb{R}^1. $$ The identity
(\ref{5.3}) with the fixed $\eta$ we rewrite in the form
\begin{align}\label{5.11}
\underset{Q^+_{2R}(z^0)}{\int}-v\,\xi_t\,dz+\underset{\Gamma_{2R}(z^0)}{\int}-v^G\,\xi_t\,d\Gamma=
\underset{Q^+_{2R}(z^0)}\int[(\Phi,\nabla\,\xi)+F\,\xi]\,dz+\underset{\Gamma_{2R}(z^0)}
{\int}[(\Phi^G,\nabla'\xi)+F^G\,\xi]\,d\Gamma
\end{align}
where
\begin{align}\label{5.12}
\Phi(z)=(a(z)\nabla\,u)\,\omega(t)\,d(x),\,\,\,\,\,\Phi^G(z')=(b(z')\nabla'u^G)\,\omega(t)\,d(x',0),\\
\notag F(z)=(u(z)-u^G_{R,z^0})\,\omega'(t)\,d(x)
+(a(z)\nabla\,u,\nabla\,d)\,\omega(t)+f(z)\omega(t)\,d(x), \\
\notag
  F^G(z')=(u^G-u^G_{R,z^0})\omega'(t)\,d(x',0)
 +(b(z')\nabla'u^G,\nabla'd(x',0))\omega(t)+\psi(z')\,\omega(t)\,d(x',0).
 \end{align}
We put $f(z)$ and $\psi(z')=0$  for $t\in \mathbb{R}^1\setminus
\Lambda_{2R}(t^0)$  and remark that  the functions
$v,\,v^G,\,\Phi,\,\Phi^G,$  $F,$ and $F^G$  vanish for $t\in
\mathbb{R}^1\setminus \Lambda_{2R}(t^0)$.  The identity
(\ref{5.11}) can be written in the form
\begin{align}\label{5.13}
\int_{\hat{Q}}-v\,\xi_t\,dz+\int_{\hat{\Gamma}}-v^G\,\xi_t\,d\Gamma
=\int_{\hat{Q}}[(\Phi,\nabla\,\xi)+F\,\xi]\,dz
\\ \notag + \int_{\hat{\Gamma}}[(\Phi^G,\nabla'\xi)+
F^G\,\xi]\,d\Gamma,\,\,\,\forall \xi\in W^{1,1}_2(\hat{Q})\cap
W^{1,1}_2(\hat{\Gamma}).
\end{align}
For any $w(t)\in L^1(\mathbb{R}^1)$ we define the Steklov averages
$$w_h(t)=\frac{1}{h}\int_t^{t+h}w(\tau)\,d\tau,\quad
w_{\overline{h}}(t)=\frac{1}{h}\int_{t-h}^t w(\tau)\,d\tau.
$$
We put $\xi(z)=g_{\overline{h}}(x,t)$ in (\ref{5.13})  for any
$g\in W^{1,1}_2(\hat{Q})\cap W^{1,1}_2(\hat{\Gamma})$.

It allows us to transform  (\ref{5.13})  in the way:

\begin{align}\label{5.14}
\int_{\hat{Q}}-v_h\,g_t\,dz+\int_{\hat{\Gamma}}-v^G_h\,g_t\,d\Gamma=\int_{\hat{Q}}[(\Phi_h,\nabla\,g)+
F_h\,g]\,dz +
\int_{\hat{\Gamma}}[(\Phi^G_h,\nabla'g)+F^G_h\,g]\,d\Gamma.\end{align}
If to fix $g(z)=\chi(t)\,\theta(x)$,  $\chi\in
C^\infty_0(\mathbb{R}^1)$  and $\theta\in W^1_2(B^+_{2R}(x^0))\cap
W^1_2(\gamma_{2R}(x^0)),\,\,\theta|_{S^+_{2R}(x^0)}=0$, then
\begin{align}
\label{5.15}
-\int_{\mathbb{R}^1}\chi'(t)[(v_h,\theta)_{2,B^+_{2R}}+(v_h^G,\theta)_{2,\gamma_{2R}}]\,dt
 \\ \notag
 =\int_{\mathbb{R}^1}\chi(t)[(\Phi_h,\nabla\,\theta)_{2,B^+_{2R}}+(F_h,\theta)_{2,B^+_{2R}}
 +(\Phi^G_h,\nabla'\theta)_{2,\gamma_{2R}}+(F^G_h,\theta)_{2,\gamma_{2R}}]\,dt\end{align}
 for any $\chi\in C^\infty_0(\mathbb{R}^1).$

 By the definition of the weak derivative,
 \begin{align}\label{5.16}
 \frac{d}{dt}\left[(v_h,\theta)_{2,B^+_{2R}}+(v^G_h,\theta)_{2,\gamma_{2R}}\right]=K_h(\theta)
 \end{align}
 for almost all $t\in \mathbb{R}^1$,  here
$$K_h(\theta):=(\Phi_h,\nabla\theta)_{2,B^+_{2R}}+(F_h,\theta)_{2,B^+_{2R}}+(\Phi^G_h,\nabla'\theta)_{2,\gamma_{2R}}+(F^G_h,\theta)_{2,\gamma_{2R}}.
$$
If we fix $\theta\in \overset{0}{W}^{1}_2(B^+_{2R}(x^0))$ in
(\ref{5.15})  then
\begin{align}\label{5.16'}
\frac{d}{dt}(v_h,\theta)_{2,B^+_{2R}}=(\Phi_h,\nabla\theta)_{2,B^+_{2R}}+(F_h,\theta)_{2,B^+_{2R}}
\end{align}
for almost all $t\in \mathbb{R}^1.$

Now using  well-known properties of the Steklov averages  (see,
for example, \cite{LSU}, Ch.2, Lemma 4.7, Ch.3, Lemma 4.1)   we
obtain the relation
\begin{align}\label{5.17}
(\frac{d}{dt}v_h,\theta)_{2,B^+_{2R}}+(\frac{d}{dt}v^G_h,\theta)_{2,\gamma_{2R}}=K_h(\theta),\,\,\,\,\,a.a.
\,\,t\in \mathbb{R}^1.
\end{align}
By $\tilde{p}(\alpha)$ we denote the Fourier transformation of a
function $p\in L^1(\mathbb{R}^1)$:
$$\tilde{p}(\alpha)=\frac{1}{\sqrt{2\pi}}\int_{\mathbb{R}^1}p(t)\exp^{-i\alpha\,t}\,dt.$$
We apply the Fourier transformation in $t$ to equality
(\ref{5.17}) and get
$$
-i\,\alpha(\widetilde{v_h(x,\alpha)},\theta)_{2,B^+_{2R}}-i\,\alpha(\widetilde{v_h^G(x,\alpha)},\theta)_{2,\gamma_{2R}}
=\widetilde{K_h(\theta)}.
$$
Multiplying the last relation by $i\,sign\,\alpha$ and putting in
it $\theta=\widetilde{v_h(x,\alpha)}$  we obtain the relation

\begin{align*}
|\alpha|\,\|\widetilde{v_h(x,\alpha)}\|^2_{2,B^+_{2R}}+|\alpha|\,\|\widetilde{v^G_h(x',\alpha)}\|^2_{2,\gamma_{2R}}
=
i\,sign\,\alpha[(\widetilde{\Phi_h},\widetilde{\nabla\,v_h})_{2,B^+_{2R}}+(\widetilde{F_h},\widetilde{v_h})_{2,B^+_{2R}}
\end{align*}
\begin{align*}
+(\widetilde{\Phi^G_h},\widetilde{\nabla'v^G_h})_{2,\gamma_{2R}}+(\widetilde{F^G_h},\widetilde{v^G_h})_{2,\gamma_{2R}}].
\end{align*}
Now we integrate the last relation in $\alpha\in \mathbb{R}^1$ and
derive the following equality:
\begin{align}\label{5.18}
J_h:=\int_{\mathbb{R}^1}|\alpha|\,\|\widetilde{v_h(x,\alpha)}\|^2_{2,B^+_{2R}}d\alpha
+\int_{\mathbb{R}^1}|\alpha|\,\|\widetilde{v^G_h(x',\alpha)}\|^2_{2,\gamma_{2R}}d\alpha
\\ \notag =\int_{\mathbb{R}^1}i\,sign\,\alpha\{(\widetilde{\Phi_h},\widetilde{\nabla\,v_h})_{2,B^+_{2R}}
+(\widetilde{F_h},\widetilde{v_h})_{2,B^+_{2R}}
+(\widetilde{\Phi^G_h},\widetilde{\nabla'v^G_h})_{2,\gamma_{2R}}+(\widetilde{F^G_h},\widetilde{v^G_h})_{2,\gamma_{2R}}\}d\alpha.
\end{align}
The Parseval's equality provides the estimate
\begin{align}\label{5.19}
J_h\leq
\|\Phi_h\|_{2,\hat{Q}}\|\nabla\,v_h\|_{2,\hat{Q}}+\|\Phi^G_h\|_{2,\hat{\Gamma}}\|\nabla'v^G_h\|_{2,\hat{\Gamma}}
\\ \notag
+R^2(\|F_h\|^2_{2,\hat{Q}}+\|F^G_h\|^2_{2,\hat{\Gamma}})+R^{-2}(\|v_h\|^2_{2,\hat{Q}}+\|v^G_h\|^2_{2,\hat{\Gamma}}).
\end{align}
Let now $h\to 0 $ in (\ref{5.19}).  Then  $$J=\lim_{h\to 0} J_h=
\|v\|^2_{H^{1/2}(\mathbb{R}^1;L^2(B^+_{2R}(x^0)))}+
\|v^G\|^2_{H^{1/2}( \mathbb{R}^1;L^2(\gamma_{2R}(x^0)))}.$$ Thus,
\begin{align}\label{5.20}
J\leq
\|\Phi\|_{2,\hat{Q}}\|\nabla\,v\|_{2,\hat{Q}}+\|\Phi^G\|_{2,\hat{\Gamma}}\|\nabla'v^G\|_{2,\hat{\Gamma}}+
R^2(\|F\|^2_{2,\hat{Q}}+\|F^G\|^2_{2,\hat{\Gamma}})+R^{-2}(\|v\|^2_{2,\hat{Q}}+\|v^G\|^2_{2,\hat{\Gamma}}).
\end{align}
Using  definitions (\ref{5.8}) and (\ref{5.12})  of the functions
$v,\,v^G,\, \Phi,\,\Phi^G,\, F\,$ and  $F^G$,  we derive estimates
(\ref{5.9}) and (\ref{5.10}). $\bullet$
\medskip

The next step is to derive the energy estimate for a weak solution
$u^0$ of problem (\ref{5.6}).

\bigskip

\noindent{\bf Proposition \,5.5.}\,\, \textit{Let the assumptions
{\bf H1--H4 }  and  {\bf H5'} hold, $u^0$ be a weak solution to
problem (\ref{5.6}).  For any fixed $ Q^+_{2R}(z^0),$ $z^0\in
\Gamma_{1-q}(0),$ $2R<q$, where  $q\in (0,1)$,  the following
estimate}
\begin{align}\label{5.21}
\int_{Q^+_R(z^0)}|\nabla\,u^0|^2\,dz\leq
\frac{c}{R^2}\int_{Q^+_{2R}(z^0)}|u^0|^2\,dz+
c\,\mathbb{L}_0\,R^{n+2\alpha}, \end{align} \textit{ is valid
where  $\mathbb{L}_0$  is defined in (\ref{E10*}).}

\medskip

\noindent {\bf Proof.} \,\, We fix a number $q\in (0,1)$,  a
cylinder $Q^+_{2R}(z^0),\,\,z^0\in \Gamma_{1-q}(0)$, $2R<q$,  and
put in (\ref{5.7})  $\eta=u^0(z)\,\omega^2(t)\,d^2(x)\in
\overset{0}{W}^{1,1/2}_2(Q^+_{2R}(z^0))$,  the functions $d(x)$
and $\omega(t)$ such as in (\ref{5.8}).  Let  $v$ be the function
defined by (\ref{5.8}),  we put
\begin{align}\label{5.22}
v^0(z)=u^0(z)\,\omega(t)\,d(x),\quad
v^G(z')=\hat{u^G}(z')\,\omega(t)\,d(x',0),\,\,\,\hat{u^G}=u^G(z')-u^G_{2R,z^0},
\\ \notag
w^G(z)=\hat{u^G}(z')\,\omega(t)\,d(x),  \quad v(z)=v^0(z)+w^G(z).
\end{align}
After trivial calculations in (\ref{5.7}) with the fixed $\eta$ we
use estimate (\ref{E10'}) and  obtain the inequalities
\begin{align}\label{5.23}
\int_{Q^+_{2R}(z^0)}|\nabla\,v^0|^2\,dz\leq
|J_{2R}|+\frac{c}{R^2}\int_{Q^+_{2R}}|u^0|^2\,dz+\frac{c}{R}\int_{\Gamma_{2R}}|u^G-u_{2R}^G|^2\,d\Gamma
 \\ \notag +c\,R\,\int_{\Gamma_{2R}}|\nabla'u^G|^2\,d\Gamma
+c\,R^2\int_{Q^+_{2R}}|f|^2\,dz \leq
|J_{2R}|+\frac{c}{R^2}\int_{Q^+_{2R}}|u^0|^2\,dz+c\,\mathbb{L}_0
\,R^{n+2\alpha}
\end{align}
where the integral
$$J_{2R}=\int_{Q^+_{2R}}u^G_t(z')\,\omega(t)d(x)v^0(z)\,dz$$ we estimate
below.

\bigskip

As $u^G\in C^\alpha(\Gamma_{1-q}(0))$,
\begin{align*}
\left|\int_{Q^+_{2R}(z^0)}\hat{u^G}(z')v^0(z)\omega'(t)\,d(x)\,dz\right|\leq
\frac{c}{R^2}\int_{Q^+_{2R}(z^0)}|u^0|^2\,dz+c\,\mathbb{L}_0\,R^{n+2\alpha},\end{align*}
and it follows from (\ref{5.23})  that
\begin{align}\label{5.24}
\int_{Q^+_{2R}(z^0)}|\nabla\,v^0|^2\,dz\leq
c\,|I_{2R}|+\frac{c}{R^2}\int_{Q^+_{2R}(z^0)}|u^0|^2\,dz+c\,\mathbb{L}_0\,R^{n+2\alpha},
\end{align}
where
$$|I_{2R}|=\left|\int_{Q^+_{2R}(z^0)}w^G_t(z)\,v^0(z)\,dz\right|.$$

To estimate $|I_{2R}|$  we go back to the proof of the Proposition
5.4 and apply the Fourier transformation (with respect to the
variable $t$) to the relation (\ref{5.16'}). Then
$$
-i\,\alpha(\widetilde{v_h(x,\alpha)},
\theta(x))_{2,B^+_{2R}(z^0)}=(\widetilde{\Phi_h(x,\alpha)},\nabla\,\theta)_{2,B^+_{2R}}
+(\widetilde{F_h(x,\alpha)},\theta)_{2,B^+_{2R}},
$$
for any $\theta\in \overset{0}{W}^1_2(B^+_{2R}(x^0)).$

We multiply the last relation by $i\,sign\,\alpha $  and put
$\theta=\widetilde{v^0_h(x,\alpha)}$. It follows  that
\begin{align*}
|\alpha|\,(\widetilde{v_h(x,\alpha)},\widetilde{v_h^0(x,\alpha)})_{2,B^+_{2R}}=
i\,sign\,\alpha\,
[(\widetilde{\Phi_h},\widetilde{\nabla\,v^0_h})_{2,B^+_{2R}}+(\widetilde{F_h},\widetilde{v^0_h})_{2,B^+_{2R}}]
\end{align*}
As $\widetilde{v_h}=\widetilde{v^0_h}+\widetilde{w^G_h}$, we
obtain (after integrating in $\alpha\in \mathbb{R}^1$  the last
relation) that
\begin{align}\label{5.25}
\int_{\mathbb{R}^1}|\alpha|\,\|\widetilde{v^0_h(x,\alpha)}\|_{2,B^+_{2R}}^2\,d\alpha
+\int_{\mathbb{R}^1}|\alpha|\,(\widetilde{w^G_h},\widetilde{v^0_h})_{2,B^+_{2R}}\,d\alpha\leq
\|\widetilde{\Phi_h}\|_{2,\hat{Q}}\,\|\widetilde{\nabla\,v^0_h}\|_{2,\hat{Q}}
+\|\widetilde{F_h}\|_{2,\hat{Q}}\,\|\widetilde{v^0_h}\|_{2,\hat{Q}}.
\end{align}
Taking into account the inequality
\begin{align*}
\left|\int_{\mathbb{R}^1}|\alpha|\,(\widetilde{w_h^G},\widetilde{v^0_h})_{2,B^+_{2R}}\,d\alpha\right|\leq
\frac{1}{2}\int_{\mathbb{R}^1}|\alpha|\,\|\widetilde
{w^G_h}\|^2_{2,B^+_{2R}}\,d\alpha+\frac{1}{2}\int_{\mathbb{R}^1}
|\alpha|\,\|\widetilde{v^0_h}\|^2_{2,B^+_{2R}}\,d\alpha,
\end{align*}
we derive  from (\ref{5.25})  that
$$\int_{\mathbb{R}^1}|\alpha|\,\|\widetilde{v^0_h}\|_{2,B^+_{2R}}^2\,d\alpha \leq\int_{\mathbb{R}^1}|\alpha|\,\|\widetilde{w^G_h}\|^2_{2,B^+_{2R}}\,d\alpha+
2\,\|\widetilde{\Phi_h}\|_{2,\hat{Q}}\|\widetilde{\nabla\,v^0_h}\|_{2,\hat{Q}}
+2\,\|\widetilde{F_h}\|_{2,\hat{Q}}\,\|\widetilde{v^0_h}\|_{2,\hat{Q}}.
$$
If $h\to 0$ in the last inequality then
\begin{align}\label{5.26}
\|v^0\|^2_{H^{1/2}(\mathbb{R}^1,L^2(B^+_{2R}(x^0)))}\leq
2\,\|\Phi\|_{2,Q^+_{2R}(z^0)}\|\nabla\,v^0\|_{2,Q^+_{2R}}+2\,\|F\|_{2,Q^+_{2R}}\|v^0\|_{2,Q^+_{2R}}
\\ \notag +\|w^G\|^2_{H^{1/2}(\mathbb{R}^1,L^2(B^+_{2R}(x^0)))}.
\end{align}
By the definition of the function $w^G$ and due to the H\"{o}lder
continuity of the function $u^G(z')$ along $\Gamma_{1-q}(0)$, we
obtain from (\ref{5.26}) that
\begin{align}\label{5.27}
\|w^G\|^2_{H^{1/2}(\Lambda_{2R}(t^0),L^2(B^+_{2R}(x^0)))}=\int_{\Lambda_{2R}(t^0)}\int_{\Lambda_{2R}(t^0)}\int_{B^+_{2R}}
d^2(x)\frac{|\hat{u}^G(x',t)\omega(t)-\hat{u}^G(x',\tau)\omega(\tau)|^2}{|t-\tau|^2}\,dx\,dt\,d\tau
\\ \notag \leq c\,\mathbb{L}_0\,R^{n+2\alpha}.
\end{align}
Using the definition (\ref{5.12}) of the functions $\Phi$ and $F$,
we derive from (\ref{5.26}) and (\ref{5.27}) the inequality
\begin{align}\label{5.28}
\|v^0\|^2_{H^{1/2}(\Lambda_{2R}(t^0),L^2(B^+_{2R}))}\leq
c\|\nabla\,v^0\|^2_{2,Q^+_{2R}}+\frac{c\,\|u^0\|^2_{2,Q^+_{2R}}}{R^2}+c\,\mathbb{L}_0\,R^{n+2\alpha}.
\end{align}
Now estimates (\ref{5.27}) and (\ref{5.28}) allows us to estimate
the expression $|I_{2R}|$:
$$
|I_{2R}|\leq\|w^G\|_{H^{1/2}(\Lambda_{2R}(t^0),L^2(B^+_{2R}))}\,\|v^0\|_{H^{1/2}(\Lambda_{2R}(t^0),L^2(B^+_{2R}))}
$$
$$\leq
\varepsilon\,\|v^0\|^2_{H^{1/2}(\Lambda_{2R}(t^0),L^2(B^+_{2R}))}+
\varepsilon^{-1}\|w^G\|^2_{H^{1/2}(\Lambda_{2R}(t^0),L^2(B^+_{2R}))}
$$
$$\leq \varepsilon\|\nabla\,v^0\|^2_{2,Q^+_{2R}}+c\,
\mathbb{L}_0(\varepsilon^{-1}+1)R^{n+2\alpha}+
\frac{c\,\|u^0\|^2_{2,Q^+_ {2R}}}{R^2},\quad \forall
\varepsilon>0.
$$
We put $\varepsilon=1/2$ in the last inequality and apply it to
estimate the right hand side of (\ref{5.24}). Thus,
\begin{align}\label{5.29}
\int_{Q^+_{2R}}|\nabla\,v^0|^2\,dz\leq
c\,\mathbb{L}_0\,R^{n+2\alpha}+ c\,R^{-2}\int_{Q^
+_{2R}}|u^0|^2\,dz.
\end{align}
As $v^0=u^0$ in $Q^+_R(z^0)$, estimate (\ref{5.21}) follows from
(\ref{5.29}). $\bullet$
\bigskip

\noindent{\bf Proof of Theorem 2.2.} \,\, Let now  $z^0\in
\Gamma_{1-q}(0)$ where the number $q\in (0,1)$ is fixed
arbitrarily, and a cylinder $Q^+_{2R}(z^0)\subset Q^+_1(0)$.

We put
 $$\Psi(\rho,z^0)=\vint_{Q_\rho(z^0)}|u^0(z)|^2\,dz,\quad \rho\leq
 2R,\quad
 u^0(z)=u(z)-u^G(z').
 $$
For a fixed $\varepsilon>0$ and the matrix
$$A(t)=\vint_{B^+_{R}(x^0)}a(x,t)\,dx$$  we apply Lemma 3.4  to the
function $u^0\in W^{1,0}_2(Q^+_R(z^0),\,\,u^0|_{\Gamma_R(z^0)}=0$.
By the lemma, there exist a constant $C_\varepsilon$, an
$A(t)$-caloric function $h\in
W^{1,0}_2(Q^+_{R/2}(z^0)),\,\,h|_{\Gamma_{R/2}(z^0)}=0$,  and a
function $\phi_0\in C^1_0(Q_R(z^0)), \,$
$\sup_{Q_R(z^0)}|\nabla\,\phi_0|\leq 1$,\,( $|\phi_0|\leq
c\,R,$\,$\,|(\phi _0)_t|\leq \frac{c}{R}$)  such that
\begin{align}\label{5.30}
\vint_{Q^+_{R/2}(z^0)}(|h|^2+R^2|\nabla\,h|^2)\,dz\leq
2^{n+2}\vint_{Q^+_R(z^0)}(|u^0|^2+R^2\,|\nabla\,u^0|^2)\,dz,
\end{align}
\begin{align}\label{5.31}
\vint_{Q^+_{R/2}(z^0)}|u^0-h|^2\,dz\leq
\varepsilon\vint_{Q^+_R(z^0)}(|u^0|^2+R^2\,|\nabla\,u^0|^2)\,dz+C_\varepsilon\,R^2\mathcal{L}^2_a(R,\phi_0),
\end{align}
where
$$\mathcal{L}^2_a(R,\phi_0)=\left|\vint_{Q^+_R(z^0)}(-u^0\,(\phi_0)_t+(A(t)\nabla\,u^0,\nabla\,\phi_0))\,dz\right|^2.
$$
We estimate the function $\Psi(\rho)=\Psi(\rho,z^0),\,\,\rho\leq
R/2,$  with the help of the Campanato estimate (\ref{Cam1}) for
the $A(t)$-caloric function $h$, inequalities (\ref{5.30}),
(\ref{5.31}),  and the Friedrichs inequality:
\begin{align}\label{5.32}
\Psi(\rho)\leq
2\vint_{Q^+_\rho(z^0)}|u^0-h|^2\,dz+2\vint_{Q^+_\rho(z^0)}|h|^2\,dz\leq
c\left(\frac{R}{\rho}\right)^{n+2}\vint_{Q^+_{R/2}(z^0)}|u^0-h|^2\,dz
\\ \notag +c\,\left(\frac{\rho}{R}\right)^2\vint_{Q^+_{R/2}(z^0)}|h|^2\,dz\leq
c\,[\varepsilon\,R^2\,\vint_{Q^+_{R}(z^0)}|\nabla\,u^0|^2\,dz
  +C_\varepsilon\,R^2\mathcal{L}^2_a(R,\phi_0)]\left(\frac{R}{\rho}\right)^{n+2}  \\ \notag +
 c\,\left(\frac{\rho}{R}\right)^2\,R^2\,\vint_{Q^+_R(z^0)}|\nabla\,u^0|^2\,dz.
\end{align}
By estimate (\ref{5.21}),
\begin{align}\label{5.33}
R^2\vint_{Q^+_R(z^0)}|\nabla\,u^0|^2\,dz\leq
c\,\vint_{Q^+_{2R}(z^0)}|u^0|^2\,dz+c\,\mathbb{L}_0\,R^{2\alpha}.
\end{align}
It follows from (\ref{5.32}), (\ref{5.33}) that
\begin{align}\label{5.34}
\Psi(\rho)\leq
c\,\left[\left(\frac{\rho}{R}\right)^2+\varepsilon\left(\frac{R}{\rho}\right)^{n+2}\right]\Psi(2R)
+c\,\mathbb{L}_0\,R^{2\alpha}+C_\varepsilon\,R^2\,\mathcal{L}^2_a(R,\phi_0)\left(\frac{R}{\rho}\right)^{n+2}.
\end{align}
Now we address to identity (\ref{5.7}) to estimate
$\mathcal{L}^2_a(R,\phi_0)$.  We put $$\Delta \,a=a(x,t)-A(t)$$
and attract  conditions  {\bf H2,\,H5'}  to derive the
inequalities:
$$\mathcal{L}^2_a(R,\phi_0)=\left|\vint_{Q^+_R}(\Delta\,a\,\nabla\,u^0,\nabla\,\phi_0)\,dz
+\vint_{Q^+_R}[\hat{u}^G\,(\phi_0)_t-(a\,\nabla'u^G,\nabla\,\phi_0)+f\,\phi_0]\,dz\right|^2
$$
$$\leq
\vint_{Q^+_R}|\Delta\,a|^2\,dz\vint_{Q^+_R}|\nabla\,u^0|^2\,dz
+\vint_{Q^+_R}|\hat{u}^G|^2\,dz\vint_{Q^+_R}|(\phi_0)_t|^2\,dz+c\,\vint_{\Gamma_{2R}}|\nabla'u^G|^2\,d\Gamma
$$
$$+c\,R^2\vint_{Q^+_{2R}}|f|^2\,dz\leq
q^2_a(R)\vint_{Q^+_R}|\nabla\,u^0|^2\,dz+c\,\mathbb{L}_0\,R^{-2+2\alpha}.
$$
Certainly, we have used estimate (\ref{E10'}) from Theorem 2.1.
Now it follows from  (\ref{5.34}) and (\ref{5.21}) that
\begin{align}\label{5.35}
\Psi(\rho)\leq
c\left[\left(\frac{\rho}{R}\right)^2+\varepsilon\left(\frac{R}{\rho}\right)^{n+2}
+C_\varepsilon\,q_a^2(R)\left(\frac{R}{\rho}\right)^{n+2}\right]\Psi(2R)+c\,\mathbb{L}_0\,R^{2\alpha}.
\end{align}
We put $r=2R$ in the last inequality and obtain the relation
\begin{align}\label{5.36}
\Psi(\rho)\leq
c_0\left[\left(\frac{\rho}{r}\right)^2+\varepsilon\left(\frac{r}{\rho}\right)^{n+2}
+C_\varepsilon\left(\frac{r}{\rho}\right)^{n+2}\,q_a^2(r)\right]\Psi(r)+c_1\,r^{2\alpha}\mathbb{L}_0,\quad\forall
\, \rho\leq r/4.\end{align}

Now we fix  $\beta=\frac{1+\alpha}{2}>\alpha $  and put
$\rho=\tau\,r,\,\,\tau\leq 1/4$ in (\ref{5.36}):

\begin{align}\label{5.37}
\Psi(\tau\,r)\leq
c_0[\tau^2+\varepsilon\tau^{-(n+2)}+C_\varepsilon\tau^{-(n+2)}q_a^2(r)]\Psi(r)+c_1\,\mathbb{L}_0\,r^{2\alpha}.
\end{align}
Further we fix $\tau\leq 1/4$ such that
\begin{align}\label{5.38}
c_0\,\tau^2<\frac{\tau^{2\beta}}{4}.
\end{align}
Then we put  $\varepsilon>0$ to satisfy the inequality
\begin{align}\label{5.39}
c_0\,\varepsilon\,\tau^{-(n+2)}<\frac{\tau^{2\beta}}{4}.
\end{align}
At last, we fix $r_0$:
\begin{align}\label{5.40}
c_0\,C_\varepsilon\,\tau^{(n+2)}q^2(r_0)<\frac{\tau^{2\beta}}{4}.
\end{align}
Under conditions  (\ref{5.38}) - (\ref{5.40}) we have the
inequality
\begin{align}\label{5.41}
\Psi(\tau\,r)<\tau^{2\beta}\Psi(r)+c_1r^{2\alpha}\mathbb{L}_0.
\end{align}

For the fixed $\tau,\varepsilon,\,r_0$, we can change $r$ by
$\tau^j r, \,\,j\in \mathbb{N},$
 and repeat all considerations. In a result,
  we get
\begin{align}\label{5.42}
\Psi(\tau^{j+1}r)\leq
\tau^{2\beta}\Psi(\tau^jr)+c_1\,\tau^{2\alpha\,j}\mathbb{L}_0\,r^{2\alpha}.
\end{align}
The iterating process provides that
\begin{align}\label{5.43}
\Psi(\tau^{j+1}r)\leq c_1\,\tau^{2\beta
j}\Psi(r)+c_2\,\tau^{2\alpha j}\,\mathbb{L}_0\,r^{2\alpha}.
\end{align}
The constants $c_1$ and $c_2$ in (\ref{5.42}) and (\ref{5.43}) do
not depend on $z^0$  and  $r$.

It follows from (\ref{5.43}) that
\begin{align}\label{5.44}
\Psi(\rho)\leq
c\rho^{2\alpha}\left(\frac{\Psi(r)}{r^{2\alpha}}+\mathbb{L}_0\right),\quad
\forall\rho\leq r\leq r_0.
\end{align}
Thus, for any $z^0\in\Gamma_{1-q}(0)$
\begin{align}\label{5.45}
\sup_{\rho\leq
r_0}\frac{1}{\rho^{n+2+2\alpha}}\int_{Q^+_\rho(z^0)}|u^0|^2\,dz\leq
c\{r_0^{-(n+2+2\alpha)}\int_{Q^+_{r_0}(z^0)}|u^0|^2\,dz+\mathbb{L}_0\}.
\end{align}
Taking into account that $u^G\in
C^\alpha(\Gamma_{1-q}(0);\delta)$, we derive from (\ref{5.45})
that
$$\Phi(\rho,z^0):=\frac{1}{\rho^{n+2+2\alpha}}\int_{Q^+_\rho(z^0)}|u-u_{\rho,z^0}|^2\,dz
\leq\frac{1}{\rho^{n+2+2\alpha}}\int_{Q^+_\rho(z^0)}|u-u^G_{\rho,z^0}|^2\,dz
$$
$$\leq
\frac{2}{\rho^{n+2+2\alpha}}\int_{Q^+_\rho(z^0)}|u^0|^2\,dz
+\frac{2\rho}{\rho^{n+2+2\alpha}}\int_{\Gamma_\rho(z^0)}|u^G-u^G_{\rho,z^0}|^2\,d\Gamma
$$
$$ \leq
 c(r_0^{-1})\{\|u\|^2_{W^{1,0}_2(Q^+_1)}+\|u^G\|^2_{W^{1,0}_2(\Gamma_1)}\}
+c\,(\|f\|^2_{L^{2,n-2+2\alpha}(Q^+_1;\delta)}+\|\psi\|^2_{L^{2,n-3+2\alpha}(\Gamma_1;\delta)}).
$$
It follows that
\begin{align}\label{5.46}
\sup_{z^0\in \Gamma_{1-q},\,\rho\leq r_0}\Phi(\rho,z^0)\leq
\,c(r_0^{-1})\{\|u\|^2_{W^{1,0}_2(Q^+_1)}+\|u^G\|^2_{W^{1,0}_2(\Gamma_1)}\}
+c\,(\|f\|^2_{L^{2,n-2+2\alpha}(Q^+_1;\delta)}
\\ \notag +\|\psi\|^2_{L^{2,n-3+2\alpha}(\Gamma_1(0);\delta)}).
\end{align}
We recall that in estimate (\ref{5.46}) the number $r_0\leq q$
depends on the data only.
\medskip

 At the same time,  it follows from Proposition 5.2 that
 $\Phi(\rho,\xi)$ is estimated  for $\xi\in
 Q^+_{1-q}(0)$ and $\rho<\delta(\xi,\Gamma_1)$ (see (\ref{5.4})). Usual "sewing "
 procedure allows us to derive the estimate
 \begin{align}\label{5.47}
 \sup_{\rho\leq r_0,z^0\in
 \overline{Q^+_{1-q}(0)}}\Phi(\rho,z^0)\leq c\,\,\mathbb{L}_0.
 \end{align}
  It means that the seminorm of $u$ in
 $\mathcal{L}^{2,n+2+2\alpha}(Q^+_{1-q}(0);\delta)$ is estimated.
 By the isomorphism of this space to
 $C^\alpha(\overline{Q^+_{1-q}(0)};\delta)$  we have got
 estimate of the H\"{o}lder norm of $u$. Further, using estimate (\ref{5.21}), we derive
 (\ref{E10*}). $\bullet$

 \bigskip

 \textbf{Proof of Theorem 2.3}\,\, To prove the assertion {\bf I}
 we should repeat  proof of  Theorem 2.1 up to  the relation
 (\ref{4.14}).
 In this case
 \begin{align}\label{5.48}
 K(r)\leq K_0\,r^{2\alpha_0},\quad\,
 K_0=\|f\|^2_{L^{2,\lambda_0}(Q^+_1;\delta)}+\|\psi\|^2_{L^{2,\lambda_0}(\Gamma_1;\delta)},
 \end{align}
 and inequality (\ref{4.22}) is valid  with $\alpha_0=\frac{\lambda_0+1}{2}$ and $K_0$ defined  by
 (\ref{5.48}).
 Indeed, \,taking into account the definition (\ref{kk}) of the
 expression  $K(r)$  and the assumptions   on $f$ and $\psi$,  we
 obtain validity of (\ref{5.48}) and  the assertions i) and ii) of Theorem 2.3.

 To prove the assertion {\bf II} of Theorem 2.3   we repeat all steps of the proof of Theorem 2.2 with $f\in
L^{2,2\alpha_0}(Q^+_1;\delta)$  and $\psi\in
L^{2,\lambda_0}(\Gamma_1;\delta)$ where  $\lambda_0\in (0,1)$ and
$\alpha_0=\frac{\lambda_0+1}{2}$. $\bullet$
\bigskip

The author was supported  by RFFI, grant 15-01-07650, and by grant
SPbGU, no. 6.38.670.2013.

\bigskip

\bigskip

\textit{St-Petersburg State University, Math. and Mech. Faculty,
198504, University Street, 28, Petergof, St. Petersburg,  Russia;
 e-mail: arinaark@gmail.com}

\end{document}